\documentclass[11pt]{article}
\usepackage{amsmath}
\usepackage{amssymb}
\usepackage{amsthm}
\usepackage{enumerate}
\usepackage{graphicx}

\date{\today}

\def\n{{\bf n}}
\def\R{\mathbb{R}}
\def\N{\mathbb{N}}
\def\D{\partial}

\def\eps{\varepsilon}

\newcommand{\NN}{{\mathbb N}}

\newcommand{\RR}{{\mathbb R}}

\newcommand\cC{\mathcal{C}}

\newcommand\cF{\mathcal{F}}
\newcommand\cH{\mathcal{H}}

\newcommand\cL{\mathcal{L}}
\newcommand\cM{\mathcal{M}}

\newcommand\cS{\mathcal{S}}
\newcommand\cT{\mathcal{T}}

\newcommand\cV{\mathcal{V}}

\newcommand\cZ{\mathcal{Z}}

\def\bb{{\bf b}}

\def\bn{{\bf n}}
\def\bu{{\bf u}}

\def\bv{{\bf v}}

\def\bF{{\bf F}}

\def\b1{{\bf 1}}

\def\ta{\tilde{a}}

\def\tOmega{\tilde{\Omega}}
\def\tGamma{\tilde{\Gamma}}
\def\tSigma{\tilde{\Sigma}}

\def\ov{\overline{v}}
\def\div{\mbox{{\rm div}}\;}
\def\rot{\mbox{{\rm curl}}\;}

\def\Dn{\partial _{{\bf n}}}

\newtheorem{theorem}{Theorem}[section]
\newtheorem{proposition}{Proposition}[section]
\newtheorem{lemma}{Lemma}[section]

\newtheorem{remark}{Remark}[section]

\newtheorem{stepbd}{Step}

\numberwithin{equation}{section}

\begin{document}
\begin{center}
{ \Large \bf On the ferromagnetism equations
with large variations solutions  }
\end{center}

\begin{center}

Olivier  Gu\`es and Franck Sueur\footnote{Laboratoire d'Analyse,
de Topologie et de Probabilit\'e. Centre de Math\'ematiques et
d'Informatique. 39, rue F. Joliot Curie  13453 Marseille Cedex 13 
\\ gues@cmi.univ-mrs.fr  , fsueur@cmi.univ-mrs.fr 
\\ $2000$ MSC: $35$K$50$,  $35$K$55$,   $35$Q$60$,  $35$Q$20$ 
\\  Key words: Landau-Lifschitz
equations, BKW method, Asymptotic expansion, Boundary layer}
 
\end{center}

\medskip

\bigskip


\begin{abstract}
We exhibit some large variations solutions of the Landau-Lifschitz
equations as the exchange coefficient $\eps^2$ tends to zero. 
These solutions  are
described by some asymptotic expansions which involve
some internals layers by means  of some large amplitude fluctuations  in
a neighborhood of  width  $\sim \eps$ of an hypersurface contained in the domain.
Despite the nonlinear  behaviour of these layers  we manage to justify
locally in time these asymptotic expansions.
\end{abstract}

\section{Introduction}

 Ferromagnetic materials can attain a large magnetization under the action of a small applied magnetic field.
 To explain this phenomenon, in $1907$, Weiss suggested that a \textit{spontaneous magnetization} occurs.
 In $1928$ Heisenberg explained the spontaneous magnetization postulated by Weiss in terms of the \textit{exchange energy}.
  In $1935$ Landau and Lifschitz (cf. \cite{LL}) proposed a quantitative theory, now known as \textit{micromagnetics}. 
For a  piece of ferromagnet-which is supposed to be a regular bounded open set $\Omega$ in $\RR^3$ 
with a smooth boundary, and locally on one side of $\Gamma$-
  the magnetic state at a point $x \in \Omega$ at time  $t$ is given by a
  vector $u(t,x) \in \RR^3$ which belongs to the unit sphere of $\RR^3$,
   called the  \textit{magnetic moment}.  The Landau-Lifschitz equations  read:
\begin{eqnarray}  
\label{LL1} 
 \D_t{u^\eps } =   u^\eps  \wedge ( \cH(u^\eps ) + \eps^2 \Delta u^\eps  ) 
- u^\eps  \wedge \big( u^\eps  \wedge ( \cH(u^\eps ) + \eps^2 \Delta u^\eps  ) \big) \quad  
 \mbox{ in  }  \Omega  ,
 \end{eqnarray} 
where $\eps >0$ is the exchange coefficient. 
We denote ${\cH}(u) := H_{|\Omega }  \in L^2(\Omega;\RR^3 )$ 
where the magnetic field $ H \in L^2(\RR^3;\RR^3)$, 
is the unique solution of the following elliptic problem,

\begin{equation} 
\label{H} 
\left\{\begin{array}{l} 
H \in L^2( \RR^3;\RR^3) \ ,\\ 
\\ 
\rot H =0\mbox { in }\RR^3 \ ,\\ 
\\ 
\div (H +\overline{u})=0\mbox { in }\RR^3 \ , 
\end{array} \right. 
\end{equation} 

where $\overline{u}$ means the extension of $u$ 
by $0$ outside of the set $\Omega$. 
The equations (\ref{LL1}) are supplemented by the homogeneous Neumann boundary condition:
\begin{eqnarray} 
 \label{LL2} 
\Dn u^\eps  = 0 \,  \mbox{ in \  }  \Gamma   , 
\end{eqnarray} 
where $\bn$ is the unitary outward normal at the boundary $\Gamma$, and by an initial condition:
\begin{eqnarray} 
\label{LL3} 
u^\eps _{| t=0}  = u_0 .
\end{eqnarray} 
The solution 
must also satisfy the constraint
\begin{equation} \label{unit}
| u^\eps (t,x) | = 1 , \quad  \mbox{ for \  }  x\in \Omega , \ t \geq 0 
\end{equation}
which is obviously propagated from the initial data
as soon as it is satisfied at $t=0$.

In this paper we study the asymptotic behaviour of the solutions of the Landau-Lifschitz 
equations (\ref{LL1})-(\ref{LL2})-(\ref{LL3}) as the exchange coefficient
$\eps$ tends to zero. 
From a formal point of view, when $\eps=0$, the system (\ref{LL1})-(\ref{LL2})-(\ref{LL3})
becomes 
\begin{equation}  
\label{LL0} 
\left\{ 
\begin{array}{l}  
 \D_t{u^0} =   u^0 \wedge \cH(u^0) 
- u^0 \wedge ( u^0 \wedge  \cH(u^0)  ) \quad  
 \mbox{ in  } \Omega  \\ 
\\ 
u^0_{| t=0}  = u_0 \ ,
\end{array} 
\right. 
\end{equation} 
where no boundary condition is needed. 
In the paper \cite{CFG1} it is proved that, for 
\emph{smooth enough solutions} the system (\ref{LL0})
is a "good approximation" of the full system (\ref{LL1})-(\ref{LL2})-(\ref{LL3})
in the sense that the solution $u^0$
of (\ref{LL0}) is indeed limit in $L^2([0,T]\times \Omega)$
 of solutions $u^\eps$ of (\ref{LL1})-(\ref{LL2})-(\ref{LL3})
as $\eps \rightarrow 0$. However, this result holds under
the assumption that   $u^0$ belongs
to the space $\cC\big( [0,T],H^5(\Omega)\big)$ where $H^5(\Omega)$
is the usual Sobolev space.
In particular this assumption excludes the case where $u^0$
is \emph{discontinuous } across 
an hypersurface contained in $\Omega$
and it was one  motivation  behind this paper to treat that case.

First of all, let us observe that the system (\ref{LL0})
actually admits discontinuous solutions. To simplify,
we will restrict the analysis to piecewise smooth solutions. 
We assume that   $\Sigma$ is  a smooth
compact hypersurface contained in $\Omega$.
For $0 \leq   s < \infty$ call $H^s_\Sigma(\Omega)$ the set of functions
$u\in L^2(\Omega)$ such that $u_{|\Omega_\pm} \in H^s( \Omega_\pm )$
where $H^s(\Omega_\pm)$ is the usual Sobolev space on $L^2$. 
We endow  $H^s_\Sigma(\Omega)$  with the norm
$$
\| u \|_{H^s_\Sigma} := \| u_{|\Omega_-}  \|_{H^s(\Omega_-)} 
+ \| u_{|\Omega_+}  \|_{H^s(\Omega_+)} 
$$
This definition
extends to the case when $s=\infty$: the space $H^\infty_\Sigma(\Omega)$  is
the natural Fr\'echet space. 
We get the following result of global existence of solution of
$(\ref{LL0})$ discontinuous through the hypersurface  $\Sigma$.
 \begin{theorem}
\label{hyp}
Let $s \in \,]\frac{3}{2} , \infty ]$ and $u_0 \in H^s_\Sigma(\Omega) $.
Then there  exists a unique $u^0 \in \cC^\infty\big(\RR, H^s_\Sigma(\Omega) \big)$
solution of the Cauchy problem $(\ref{LL0})$.
\end{theorem}
\begin{proof}
This result can be easily
obtained by following the proof of Proposition $4.1$ of \cite{CFG1}, with
only a few adaptations.
By the way in the closer setting of semilinear symmetric hyperbolic system, it is
well known since the works of M\'etivier \cite{m1}
that there exist some local piecewise regular solutions
discontinuous across a smooth hypersurface which is a characteristic
hypersurface of constant multiplicity for this hyperbolic system. 
In the present setting, the proof is in fact
simpler since the hypersurface $\Sigma$ is totally characteristic.
Moreover thanks to $(\ref{unit}) $ and since the operator $\mathcal{H}$ satisfies the
transmission property,  our setting allows to conclude to a global
existence. 
\end{proof}

Let us now claim a first theorem about the asymptotic behaviour of the solutions of the Landau-Lifschitz 
equations (\ref{LL1})-(\ref{LL2})-(\ref{LL3}) as the exchange coefficient
$\eps$ tends to zero.  

\begin{theorem}
\label{main}
Let $u^0 \in \cC^\infty\big(\RR, H^\infty_\Sigma(\Omega) \big) $
be a solution of $(\ref{LL0})$. 
 There exist $T>0$ and a family
of solutions $u^\eps \in \cC^\infty([0,T] \times \Omega)$,
$\eps \in ]0,1]$,
of the equation $(\ref{LL1})$ on $[0,T] \times \Omega$,
of the equation $(\ref{LL2})$ on $[0,T] \times \Gamma$,
such that  there exist $C >0$ and $\eps_0$  such that for all $\eps \in
]0,\eps_0 ]$,
\begin{eqnarray*}
|| u^\eps - u^0 ||_{L^2([0,T] \times \Omega)} \leqslant C \eps^\frac{1}{2} .
\end{eqnarray*}

\end{theorem}

To begin with some comments about Theorem $\ref{main}$ let us stress that
we do not prescribe the initial data $(\ref{LL3})$ for the $u^\eps$.
Thus the traces of the  $u^\eps$ at $t=0$ are not equal in general to  the trace of  $u^0$ at $t=0$. 
So Theorem $\ref{main}$ claims the existence of  local in time solutions  $u^\eps  \in \cC^\infty$, 
of the equation $(\ref{LL1})$ on $ \Omega$,
of the equation $(\ref{LL2})$ on $ \Gamma$, which converge to $u^0$ as
$\eps$ tends to zero in $L^2$, with a rate of convergence in $\eps^\frac{1}{2}$.

Indeed, in this paper, we will claim a more accurate result in  Theorem
$\ref{main2}$ by showing that the $u^\eps$ can be described with a WKB expansion 
which involves some boundary layers profiles. 
On one hand, a boundary layer appears near the boundary to compensate the lost of the Neumann condition
 from the complete model (\ref{LL1})-(\ref{LL2})-(\ref{LL3}) to
the limite model  (\ref{LL0})  ($\eps = 0$). 
Such a boundary layer was already studied in paper \cite{CFG1}.
 The amplitude of this boundary is
weak and its behaviour is linear.
On another hand, there are some boundary layers on each side of the hypersurface $\Sigma$. 
Their task is to compensate the lost of transmission conditions across $\Sigma$ 
from the complete model  (\ref{LL1})-(\ref{LL2})-(\ref{LL3}) to
the limit model  (\ref{LL0})  ($\eps = 0$).

\begin{remark}
\rm
Such an analysis is inspired by the paper \cite{S1} where we show that discontinuous solutions of
 multidimensional semilinear symmetric hyperbolic systems,
 which are regular outside of a smooth hypersurface characteristic of constant multiplicity, are 
  limits, when $\eps \rightarrow 0$, of solutions  $(u^\eps)_{\eps \in ]0,1]}$ of the system perturbated 
  by a viscosity of size $\eps $. 
In this paper, we adapt the method to the ferromagnetism quasi-static
model,
where in particular the non local operator ${\cH}$ occurs.
We point out that for the limit model ($\eps=0$), the hypersurface is totally characteristic.
As a consequence, the analysis involves only characteristic boundary layers.
On the opposite, \cite{S1} stresses the occurrence of characteristic and non characteristic boundary layers.
 It could be also possible -as in \cite{S1}- to study the case where the singularity  is weaker than a jump of the function $u^0$ 
as a jump of a derivative of the function $u^0$.
Then we can take $T$ as big as we want and the quality of the approximation is as better 
as the jump concerns a higher order
derivative. We also refer to papers \cite{GW}, \cite{GR}, \cite{S1} for the use of boundary layers in transmission strategy.
\end{remark}

\begin{remark}
\rm
It could be interesting to know if it is possible to obtain such a result 
in the non static case for which the Landau-Lifschitz 
equation is coupled with the Maxwell system of electromagnetic. 
For such a model an
analysis of the boundary layer induced by the Neumann boundary condition
on $\Gamma$ is performed in \cite{CFG2}.
\end{remark}

\begin{remark}
\rm
With the same method than the one used in this paper, it is
possible to get global in time
$O(\eps^s)$ approximation  for all $s < \frac{1}{2}$.
More precisely  for any $s <  \frac{1}{2}$ there exists a family
of solutions $u^\eps  \in \cC^\infty([0,T_\eps] \times \Omega)$, of the equation $(\ref{LL1})$ on $[0,T_\eps] \times \Omega$,
of the equation $(\ref{LL2})$ on $[0,T_\eps] \times \Gamma$,
$\eps \in ]0,1]$, with lim$_{\eps \rightarrow 0^+} \, T_\eps = \infty$,
such that  for all $T>0$,  there exists $C >0$ and $\eps_0$  such that for all $\eps \in
]0,\eps_0 ]$, there holds 
$|| u^\eps - u^0 ||_{L^2([0,T] \times \Omega)} \leqslant C \eps^s$ .

\end{remark}

\section{Asymptotic expansion}

Let us fix some notations.
We will use the letter $\mathcal{S}$ to denote the Schwarz space of
rapidly decreasing functions.
We  define  the boundary layer profile spaces  
\begin{eqnarray*}
  \mathcal{N}_{\pm} (T) :=    H^\infty ( [0,T] \times \Omega, \mathcal{S}( \RR^\pm ) ).
\end{eqnarray*} 
Since we will need an equation of
the boundary $\Gamma$, we fix once for all a function
$\Phi \in \cC^{\infty}(\RR^3, \RR)$ and we assume that
$\Omega = \{ \Phi > 0 \}$, $ \Gamma = \{ \Phi = 0 \}$ 
and  $| \nabla \Phi (x) | = 1 $ in an open neighborhood $\mathcal{V}_\Gamma $ 
of $\Gamma $ \footnote{ Hence for $x \in \Omega \cap \mathcal{V}_\Gamma$: 
$\Phi(x) = dist(x,\Gamma )$.}.  
Let us also  fix a function 
$\Psi \in \cC^{\infty}(\RR^3, \RR)$ such that 
$\Sigma = \{ \Psi=0 \} $
 and such that  $| \nabla \Psi (x) | = 1 $ in an open neighborhood $\mathcal{V}_\Sigma$ 
of $\Sigma $ \footnote{ Hence for $x \in \Omega \cap \mathcal{V}_\Sigma$: 
$\psi (x) = dist(x,\Sigma)$.}. 
We assume that the neighborhoods   $\mathcal{V}_\Gamma$ and  $\mathcal{V}_\Sigma$ have
been fixed small enough in order that $\cV_\Gamma \cap \cV_\Sigma =
\varnothing$.
 We will denote
$\Omega_+ := \Omega \cap \{ \Psi > 0 \}$ and
$\Omega_- := \Omega \cap \{ \Psi < 0 \}$.
We consider a $C^\infty$ unit vector field  $\partial_\bn $ which coincides on
$\mathcal{V}_\Gamma$ with  $- \nabla _x \Phi \cdot \nabla_x $ and on
$\mathcal{V}_\Sigma$ with  $- \nabla _x \Psi \cdot \nabla_x $.

\input epsf   
$$
\epsfxsize3.2in   
\epsffile{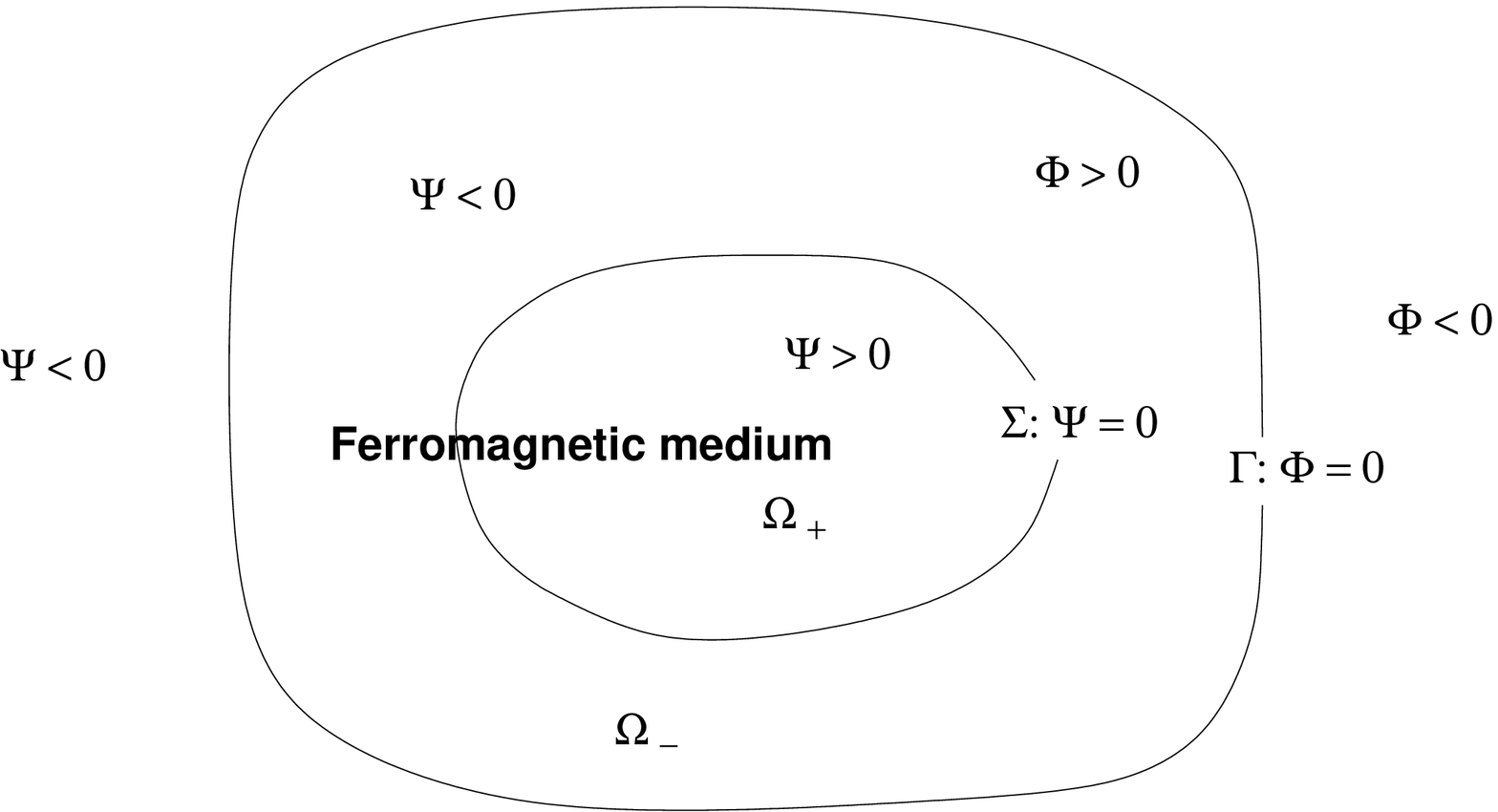}   
$$

In the easier case where $u^0$ is continuous across the hypersurface
$\Sigma$, paper \cite{CFG1} shows the existence of solutions $u^\eps $,
$\eps \in ]0,1]$,
of the equation $(\ref{LL1})$ in $\Omega$,
of the equation $(\ref{LL2})$ on $ \Gamma$,  
of the form
\begin{eqnarray*}
u^{\eps}  (t,x) &:=& u^0   (t,x) +  \eps \Big( 
\mathfrak{U} (t,x, \frac{\Phi(x)}{\eps} )  +  \mathbf{w}^{\eps}  (t,x)\Big)
\end{eqnarray*} 
where the function $\mathfrak{U}$ is in $\mathcal{N}_{+} (\infty)$ and
satisfies $\mathfrak{U}(t,x,z)=0$ for $x \notin \mathcal{V}_\Gamma$.
The function $\mathfrak{U} $ describes a boundary layer which appears near the boundary 
to compensate the lost of the Neumann condition
 from the complete model (\ref{LL1})-(\ref{LL2})-(\ref{LL3}) to
the limit model  (\ref{LL0})  ($\eps = 0$). 
 The amplitude of this boundary is
weak and its behaviour is linear.
For sake of completeness we will state this in
section $\ref{deja}$.
The  functions  $\mathbf{w}^{\eps} $ can be seen as remainders.

Here since we deal with a ground state $ u^0 $ which is  discontinuous across the  hypersurface
$\Sigma$, 
we look for solutions $u^\eps $,
$\eps \in ]0,1]$,
of the equation $(\ref{LL1})$ in $\Omega$,
of the equation $(\ref{LL2})$ on $ \Gamma$,  
of the form
\begin{eqnarray}
\label{decomp}
u^{\eps}  (t,x) :=  \mathcal{U}  (t,x, \frac{\Psi(x)}{\eps}) 
+  \eps \Big(\mathfrak{U} (t,x, \frac{\Phi(x)}{\eps} )  +  \mathbf{w}^{\eps}  (t,x)  \Big).
\end{eqnarray} 
The function $\mathcal{U}$
describes a large amplitude  internal layer
profile i.e. a sharp
transition  in the neighborhood of the hypersurface $\Sigma$ of width
$\sim \eps$.
More precisely the function $\mathcal{U}$ is $\cC^\infty$ and satisfies
\begin{eqnarray}
\label{sharp1}
\lim_{y\rightarrow \pm \infty}   \mathcal{U}  (t,x, y) &=&
  u^0  (t,x ) \quad \text{for} \ x \in \mathcal{V}_\Sigma \cap \Omega_\pm
  \\ \label{sharp2}
  \mathcal{U}  (t,x, y)&=& u^0  (t,x ) \quad \text{for} \ x \notin
  \mathcal{V}_\Sigma \ \text{and} \ y \in \R
\end{eqnarray} 
The profile $\mathfrak{U}$, as we have already said it above, was  constructed in  \cite{CFG1}. 
The  functions  $ \mathbf{w}^{\eps}  $ can still be seen as remainders.
Let us explain this time more precisely what we mean by remainders.
Let us fix a finite set of smooth vectors fields $\cT_0 = \{
\cZ_i(x;\D_x) ; i=1, \cdots ,\mu \}$ on $\RR^3$, tangent to the
surfaces $\Gamma$ and $\Sigma$ (that is satisfying $\cZ_i(x;\D_x)
\Phi =0$ on $\Gamma$ and $\cZ_i(x;\D_x)\Psi = 0$ on $\Sigma$, for
all $i\in \{1, \cdots , \mu\}$), and generating the algebra of
smooth vector fields tangent to $\Gamma \cup \Sigma$. These vector
fields can be viewed as vector fields on $\RR^{4}$ tangent to $\RR
\times \Gamma$ and to $\RR \times \Sigma$. By adding the vector
field $\D_t$ to the family, one gets the set $\cT := \{ \D_t \} \cup
\cT_0  $ which generates the set of smooth vector fields in $\RR^4 $
tangent to $(\RR \times \Gamma) \cup (\RR \times \Sigma)$. We denote
$\cZ_0 := \D_t$. For all multi-index $\alpha \in
\NN^{1+\mu}$ we note $\cZ^\alpha = \D_t^{\alpha_0}\cZ^{\alpha_1}_1.
\cdots . \cZ^{\alpha_\mu}_\mu$, with $\alpha = (\alpha_0, \alpha_1,
\cdots , \alpha_\mu)$. 
Let us introduce the usual  norm:
$$
\| u \|_{m} := \sum_{|\alpha| \leq m \, , \, \alpha \in
\NN^{1+\mu}} | \cZ^\alpha u
\|_{L^2(]0,T[ \times \Omega)},
$$
and note $H^m_{co}(]0,T[\times \Omega)$  the space of $u\in
L^2(]0,T[ \times \Omega)$ such that this norm is finite. 
We introduce the set $E$ of the family $( \mathbf{w}^\eps  )
_{0<\eps\leqslant 1})$ of functions  in $L^2(]0,T[ \times \Omega)$ such that for all $m \in \N$, there exists $\eps_0 >0$ such that
\begin{eqnarray}
\label{estiti}
\text{sup}_{0<\eps\leqslant \eps_0} \, (|| \mathbf{w}^{\eps} ||_{m} + ||\eps \partial_\n
\, \mathbf{w}^{\eps} ||_{m} +  \eps (|| \mathbf{w}^{\eps} ||_{\infty}  +
||\cZ \mathbf{w}^{\eps}
||_{\infty} + ||\eps \partial_\n \,  \mathbf{w}^{\eps}||_{\infty} ))<
\infty .
\end{eqnarray} 

In fact Theorem \ref{main} is the straightforward consequence of the
following result.

\begin{theorem}
\label{main2}

Let $u^0 \in \cC^\infty\big(\RR, H^\infty_\Sigma(\Omega) \big) $
be a solution of $(\ref{LL0})$. 
There exist $T>0$, a profile $\mathcal{U}$ in $\cC^\infty ((0,T) \times
\Omega \times \R )$ which  satisfies 
$(\ref{sharp1})-(\ref{sharp2})$
and a family $( \mathbf{w}^\eps)$ in $E$
 such that the function $u^\eps$ given by the formula (\ref{decomp}) are solutions in  $C^\infty$  
 of the equation  $(\ref{LL1})$ on $[0,T] \times \Omega$,
of the equation $(\ref{LL2})$ on $[0,T] \times \Gamma$.
\end{theorem}
Theorem \ref{main2}  exhibits some large variations solutions of the Landau-Lifschitz
equations as the exchange coefficient $\eps^2$ tends to zero, 
by means of the  asymptotic expansions (\ref{decomp}).
The remainder of the paper is devoted to the proof of Theorem \ref{main2}.
 As in \cite{CFG1}, since the magnetic moment $u$ is unimodular, the equation (\ref{LL1}) is equivalent
 for smooth solutions to the following
one:
\begin{eqnarray}
\label{LL1new}
\cL^\eps( u^\eps, \D )\, u^\eps =
    \bF\big(u^\eps, \eps \D_x u^\eps,
     \cH(u^\eps)\big)
\end{eqnarray}
where we have noted
$$
\cL^\eps(v,\D) := \D_t - \eps^2 \Delta_x - \eps^2 v \wedge \Delta_x,
$$
and
$$
\bF(u,V,H):= |V|^2u + u\wedge H - u\wedge (u\wedge H),
$$
for all $u\in \RR^3$, $V\in \cM(\RR^3,\RR^3)$, $H\in \RR^3$.
From now on we will  deal with equation (\ref{LL1new}) rather than 
(\ref{LL1}).
We will proceed in three steps.
In subsection $\ref{s1}$ we will define the profile $\mathcal{U}$  as a
local in time solution of a pair of
nonlinear equations in $\Omega \times \R_\pm $ coupled by some
transmissions conditions on $\{y=0\}$.
In  subsection $\ref{deja}$  we  will recall the results of \cite{CFG1} about the profile $\mathfrak{U}$.
In subsection $\ref{s3}$  we  will  prove the existence of some remainders $\mathbf{w}^\eps$
 till the
lifetime $T$ of the profile $\mathcal{U}$.
Eventually we will show that the remainders $\mathbf{w}^\eps$ satisfy the
uniform  estimates uniform (\ref{estiti}).

\subsection{Construction of the internal layers }
\label{s1}

Despite that $\pm \frac{\Psi(x)}{\eps}>0$ when  $x \in \Omega_\pm$ we will define $\mathcal{U}$ for all $(x,z) \in \Omega \times
\R_\pm $ since this will not cause any additional difficulty.
An Uryshon argument yields the existence of two functions  $u^{0}_\pm$  in $H^\infty ( (0,\infty) \times \Omega )$ such
that $u^{0}_\pm = u^{0}$ for all $x \in \Omega_\pm \cup  (\Omega_\mp -
\mathcal{V}_\Sigma)$.
We look for a viscous internal layer
profile $\mathcal{U}$ of the form
\begin{eqnarray}
\label{R1}
\mathcal{U}(t,x,y) := 
 \left\{
 \begin{array}{cc}
  u^{0}_+  (t,x) + \mathcal{U}_{+}  (t,x, y)  \quad
  &\text{if} \ y>0,
 \vspace{0.3cm}\\ u^{0}_-  (t,x) + \mathcal{U}_{-}  (t,x, y)  \quad
  &\text{if} \ y<0 .
 \end{array}
 \right.
\end{eqnarray} 
The functions $\mathcal{U}_\pm$ are in $\mathcal{N}_{\pm} (T)$.
These functions  describe some internal large amplitude boundary layers,
 on each side of the hypersurface $\Sigma$. 
To insure that the function $\mathcal{U}$ is in $C^1 ( (0,T) \times \Omega  \times  \R)$
it is necessary to impose the transmission conditions:
\begin{eqnarray}
 \label{p2}
\left.
\begin{array}{c}
 \mathcal{U}_{+} -  \mathcal{U}_{- } = - u^0_{+} + u^0_{- } ,
\\ \partial_y \mathcal{U}_{+} - \partial_y \mathcal{U} _{- } =0
\end{array}
\right\}
\quad \mathrm{when} \ (t,x,y) \in  (0,T) \times \Omega \times \{ 0 \}.
\end{eqnarray}
%

In Theorem $\ref{profil0}$ we will define the profiles $\mathcal{U}_\pm$ 
as local solutions of
nonlinear equations in $\Omega \times \R_\pm $ coupled by some
transmissions conditions on $\{y=0\}$.
Let us look for convenient  equations. 
We will plug the functions $u^{\eps,0}$ defined by 
$u^{\eps,0} (t,x) := \mathcal{U}(t,x, \frac{\Psi(x)}{\eps})$
instead of  $u^\eps$ in $(\ref{LL1new})$.
In general it is not possible to verify  $(\ref{LL1new}) $ but we will try
to choose the functions  $\mathcal{U}_\pm$ such that the error term is as
small as possible.
Let us begin to look at the left side of $(\ref{LL1new})$.
With (\ref{R1})
  we get  in $L^\infty$ 
\begin{eqnarray}
\label{M1}
\mathcal{L}^{\eps} (u^{\eps,0} , \partial ) u^{\eps,0}  = 
\partial_t \, u_{\pm}^0  + 
\Big( L(\mathcal{U} , \partial_t ,  \partial_y^2 )  \mathcal{U}_{\pm}   \Big)| + O(\eps) \quad   \mathrm{for} \
  x \in  \Omega_\pm  ,
\end{eqnarray}
where the vertical bar $|$ means that $y$ is evaluated in $y=
\frac{\Psi(x)}{\eps}$ and 
\begin{eqnarray*}
L(U,\partial_{t}, \partial_{y}^2) := \partial_t - \partial_y^2 - U \wedge \partial_y^2 .
\end{eqnarray*} 
We now turn to the right side of $(\ref{LL1new})$.  
We first look at  the
action of $\cH$ on the family $u^{\eps,0}$:
\begin{eqnarray*}
\cH(u^{\eps,0} ) = \cH(u^0_\pm ) - (\mathcal{U}_\pm  . n)| \, n + O(\eps).
\end{eqnarray*}
Then 
\begin{eqnarray}
\label{M2}
 \bF\big(u^\eps, \eps \D_x u^\eps,    \cH(u^\eps)\big) := 
  \bF\big(u^0_\pm, 0,   \cH(u^0_\pm)\big)
 + F_\pm ( \mathcal{U}_\pm , \partial_y  \mathcal{U}_\pm  )|  + O(\eps) \quad   \mathrm{for} \
  x \in  \Omega_\pm   ,
 \end{eqnarray}
 with for all $U\in \RR^3$, $V\in \cM(\RR^3,\RR^3)$, 
\begin{eqnarray*}
 F_\pm (  U, V  )  &:=&
 | V |^2 \, ( u^0_\pm  + U )
+ U \wedge  \cH (u^0_\pm ) - (U  . n) ( u^0_\pm  +  U)  \wedge   n
\\ \nonumber  && + U \wedge (  u^0_\pm +U)
 \wedge (   \cH (u^0_\pm ) - (U  . n) n  )
 +  u^0_\pm  \wedge (  U  \wedge  (   \cH (u^0_\pm )   - (U  . n) n ) )
\\ \nonumber  && -  (U  . n) u^0_\pm  \wedge (  u^0_\pm   \wedge    n )
\end{eqnarray*}
Thanks to (\ref{M1}) and (\ref{M2})  we get by looking at the terms at order $0$
\begin{eqnarray*}
\partial_t \,   u_{\pm}^0+ 
 L(\mathcal{U} , \partial_t ,  \partial_y^2 )  \mathcal{U}_\pm = \bF\big(u^0_\pm, 0,   \cH(u^0_\pm)\big)
 + F_\pm ( \mathcal{U}_\pm , \partial_y  \mathcal{U}_\pm  ) .
\end{eqnarray*}
Since for $x \in  \Omega_\pm$, the functions $u^0_\pm$ satisfies
$(\ref{LL0})$  we could simplify and we get
the nonlinear equations 
\begin{eqnarray}
\label{p1} L(  u_{\pm}^0 + \mathcal{U}_\pm ,\partial_{t}, \partial_{y}^2 ) \mathcal{U}_\pm =
 F_\pm ( \mathcal{U}_\pm , \partial_y  \mathcal{U}_\pm  ) .
 \end{eqnarray}
The equations  $(\ref{p1})$ are parabolic with respect to $t,y$, the
variable $x$ can be seen as a parameter.
The following theorem claims that it is possible to find some solutions
$\mathcal{U}_{\pm} \in \mathcal{N}_{\pm} (T)$ of these equations even for all $x \in   \Omega  $.

\begin{theorem}
\label{profil0}
There exists $T > 0$ and there exist some functions $\mathcal{U}_{\pm} \in \mathcal{N}_{\pm} (T)$
 which verify the equations $(\ref{p1})$ when $(t,x,y) \in  (0,T) \times \Omega \times
 \R_\pm $ and  the transmission conditions $(\ref{p2})$.
 Moreover precisely for all $x \notin \mathcal{V}_\Sigma$ and $y \in  \R_\pm $ 
there holds $\mathcal{U}_{\pm}  (t,x, y)= 0$.
\end{theorem}

\begin{proof}

We will proceed in four steps.

\begin{stepbd}
We begin to reduce the problem to homogeneous boundary conditions. 
\end{stepbd}

We introduce the functions $V_{\pm} $ and $ \mathcal{U}_{\pm} $ given by 
 the formula
\begin{eqnarray*}
V_{\pm} (t,x,y) &:=& (1-\frac{e^{\mp y}}{2}) u^0_\pm (t,x) + \frac{e^{\mp y}}{2}
u^0_\mp (t,x) ,
\\ \mathbf{W}_{\pm} (t,x,y)&:=& \mathcal{U}_{\pm} (t,x,y) \pm  \frac{1}{2} (u^0_+  (t,x)-
u^0_- (t,x)) e^{\mp y} .
\end{eqnarray*} 
Thus the transmission conditions $(\ref{p2})$ reads:
\begin{eqnarray}
 \label{np2}
\left.
\begin{array}{c}
 \mathbf{W}_{+} - \mathbf{W}_{-} =  0 ,
\\ \partial_y \mathbf{W}_{+} - \partial_y  \mathbf{W}_{-} =0
\end{array}
\right\}
\quad \mathrm{when} \ (t,x,y) \in  (0,T) \times \Omega \times \{ 0 \}.
\end{eqnarray} 
Moreover the equations  (\ref{p1})-(\ref{p2}) read for $ (t,x,y) \in 
(0,T) \times \Omega \times \R_\pm$:
\begin{eqnarray}
\label{++}
 L (V_{\pm} + \mathbf{W}_{\pm} ,\partial_{t},  \partial_{y}^2 )
\mathbf{W}_{\pm} &=& \hat{F}_{\pm} (t,x,y,\mathbf{W}_{\pm} , \partial_y
\mathbf{W}_{\pm}) ,
\end{eqnarray} 
where $\hat{F}_{\pm}$ are $C^\infty$ functions such that the functions
$\hat{F}_{\pm}(t,x,y, 0,0) $ are rapidly decreasing with respect to $y$.
\begin{stepbd}
We prove the existence of compatible initial data. 
\end{stepbd}
Let us to explain why the initial values $\mathbf{W}_{0,+}$ must satisfy some compatibility conditions 
at the corner $\{ t= y=0 \} $ are required in order to obtain smooth
solutions $\mathbf{W}_{\pm}$ of the problem  (\ref{++})-(\ref{np2}) with
$\mathbf{W}_{\pm} |_{t=0}:= \mathbf{W}_{0,\pm} $.
We start with the condition of order $0$. 
Set $t=0$ in the transmission conditions (\ref{np2}) to see that  $\mathbf{W}_{0,+}$ must satisfy  the relation 
\begin{eqnarray}
 \label{comp0}
\left.
\begin{array}{c}
 \mathbf{W}_{+} - \mathbf{W}_{0,-} =  0 ,
\\ \partial_y \mathbf{W}_{0,+} - \partial_y  \mathbf{W}_{0,-} =0
\end{array}
\right\}
\quad \mathrm{when} \ (x,y) \in  \Omega \times \{ 0 \}.
\end{eqnarray} 
Now, for each $k \geqslant 1$, apply the derivative $\partial_t^k $ to the  transmission conditions (\ref{p2}). 
We get
\begin{eqnarray*}
\left.
\begin{array}{c}
\partial_t^k \mathbf{W}_{+} - \partial_t^k\mathbf{W}_{-} =  0 ,
\\ \partial_y \partial_t^k \mathbf{W}_{+} - \partial_y \partial_t^k  \mathbf{W}_{-} =0
\end{array}
\right\}
\quad \mathrm{when} \ (t,x,y) \in  (0,T) \times \Omega \times \{ 0 \}.
\end{eqnarray*} 
Now remark that, by iteration, we can extirpate $\partial_t^k \mathbf{W}_{\pm}$ by the interior equations (\ref{p1}) in terms of
 derivatives with respect to $y$. More precisely there exists some smooth functions $C^k_{\pm} $ such that 
 $\partial_t^k \mathbf{W}_{\pm} = C^k_{\pm} ( (\partial_y^l
 \mathbf{W}_{\pm} )_{l \leqslant 2k})$.
 Thus the following $k$th order compatibility condition must hold:
 \begin{eqnarray}
 \label{compk}
\left.
\begin{array}{c}
C^k_{+} ( (\partial_y^l \mathbf{W}_{\pm} )_{l \leqslant 2k}) 
- C^k_{-} ( (\partial_y^l \mathbf{W}_{\pm} )_{l \leqslant 2k}) =0 ,
\\ \partial_y C^k_{+} ( (\partial_y^l \mathbf{W}_{\pm} )_{l \leqslant 2k}) 
- \partial_y C^k_{-} ( (\partial_y^l \mathbf{W}_{\pm} )_{l \leqslant 2k}) =0 .
\end{array}
\right\}
\quad \mathrm{when} \ (x,y) \in  \Omega \times \{ 0 \}.
\end{eqnarray} 

\begin{lemma}

There exist some initial values $\mathbf{W}_{0,\pm}$ in $H^\infty (\Omega , \mathcal{S} (\R_\pm ))$ 
which  satisfy the relation (\ref{comp0}) and (\ref{compk}) for all $k \geqslant 1$. 

\end{lemma}

\begin{proof}

As we will follows the method of \cite{S1}, we only sketch the proof for sake of completeness. 
We start by analyzing more accurately the compatibility conditions and more especially the way the functions  $C^k_{\pm} $
depend on the derivatives with respect to $y$. Indeed they exists some functions $\tilde{C}^k_{\pm} $ such that 
 \begin{eqnarray*}
C^k_{\pm} ( (\partial_y^l \mathbf{W}_{\pm} )_{l \leqslant 2k}) =   \tilde{C}^k_{\pm} ( (\partial_y^l \mathbf{W}_{\pm} )_{l \leqslant
2k-1})
+ (\partial_y^{2k} + (V_{\pm} + \mathbf{W}_{\pm} )  \wedge  \partial_y^{2k} ) \mathbf{W}_{\pm} .
\end{eqnarray*} 
Since   given two functions $ \mathbf{W}^{(0)}_{\pm}$  in $H^\infty (\Omega )$ the applications
 \begin{eqnarray*}
 \mathbf{W}_{\pm} \mapsto \mathbf{W}_{\pm} + ( V_{\pm}  +  \mathbf{W}^{(0)}_{\pm} ) \wedge   \mathbf{W}_{\pm} 
\end{eqnarray*} 
are two automorphisms of $H^\infty (\Omega )$
and an iteration, we deduce by iteration that there exists 
a family $(\mathbf{W}^{(k)}_{\pm})_{k \in \N}$ in $H^\infty (\Omega )$ such that 
 \begin{eqnarray*}
 \left.
\begin{array}{c}
C^k_{+} ( (\partial_y \mathbf{W}^{(l)}_{\pm} )_{l \leqslant 2k}) 
- C^k_{-} ( (\partial_y \mathbf{W}^{(l)}_{\pm} )_{l \leqslant 2k}) =0 ,
\\ \partial_y C^k_{+} ( (\partial_y \mathbf{W}^{(l)}_{\pm} )_{l \leqslant 2k}) 
- \partial_y C^k_{-} ( (\partial_y \mathbf{W}^{(l)}_{\pm} )_{l \leqslant 2k}) =0 .
\end{array}
\right.
 \end{eqnarray*} 
 We end the proof by a classical Borel argument.
\end{proof}

As a consequence, we will assume in the rest of the proof that the functions $\mathbf{W}_{\pm}$ vanish for $t\leqslant 0$.
\begin{stepbd}
We look for linear estimates.
\end{stepbd}
In order to use an iterative scheme, we look at the linear problem 
\begin{eqnarray}
\label{l1} 
L (\mathfrak{W}_{\pm}, \partial_{t}, \partial_{y}^2 )
\mathbf{W}_{\pm} = f_{\pm} \quad \mathrm{when} \ (t,x,y) \in  (0,T) \times \Omega \times \R_\pm ,
\\ \label{l2} 
\left.
\begin{array}{c}
 \mathbf{W}_{+} - \mathbf{W}_{-} =  0 ,
\\ \partial_y \mathbf{W}_{+} - \partial_y  \mathbf{W}_{-} =0
\end{array}
\right\}
\quad \mathrm{when} \ (t,x,y) \in  (0,T) \times \Omega \times \{ 0 \}.
\end{eqnarray} 
For all real $\lambda \geqslant 1$, the space $L^2 ((0,T) \times \Omega \times \R_\pm)$ 
is endowed with the scalar product associated
to the Euclidean norm
 \begin{eqnarray*}
 || \mathbf{W}_{\pm} ||_{0,\lambda,T} := || e^{- \lambda t} \, W_{\pm} ||_{L^2 ((0,T) \times \Omega \times \R_\pm)  }
 \end{eqnarray*} 
In order to avoid heavy notations, we will denote $W:=(\mathbf{W}_+ , \mathbf{W}_- )$, $f:=(f_+ , f_- )$ and $ \mathfrak{W} :=(
\mathfrak{W}_+ , \mathfrak{W}_- )$. 
We endow the space $ L^2 ((0,T) \times \Omega \times \R_+) \times L^2 ((0,T) \times \Omega \times \R_-)  $ with the scalar product
associated to the Euclidean norm 
 \begin{eqnarray*}
 || W ||_{0,\lambda,T} :=  || \mathbf{W}_+ ||_{+,0,\lambda,T} + || \mathbf{W}_- ||_{-,0,\lambda,T} .
 \end{eqnarray*} 
For $m \in \N$, we introduce the following weighted norms:
 \begin{eqnarray*}
 || W ||_{m,\lambda,T} := 
 \sum_{|\alpha| \leqslant m} \, || \partial_{t,x}^\alpha \, W ||_{0,\lambda,T} ,
 \ \text{and} \  | W |_{m,\lambda,T} := 
 \sum_{|\alpha| \leqslant m} \, || \partial_{t,x}^\alpha \,
 \partial_{y}^{\alpha_4} \, W ||_{0,\lambda,T},
\end{eqnarray*} 
where $\alpha:=(\alpha_0 ,..., \alpha_3) \in \N^4$ and
$\partial_{t,x}^\alpha :=\partial_{t}^{\alpha_0} \partial_{1}^{\alpha_1}
\partial_{2}^{\alpha_2} \partial_{3}^{\alpha_3}.$
\begin{proposition}

Let $R>0$. If $ \mathfrak{W}_\pm$  verify the following estimates 
 \begin{eqnarray*}
 || \mathfrak{W}_+ ||_{Lip ((0,T) \times \Omega \times \R_+)} + ||\mathfrak{W}_+ ||_{Lip ((0,T) \times \Omega \times \R_+)}
 + | \mathfrak{W} |_{m,\lambda,T} < R ,
\end{eqnarray*} 
  and the following  boundary conditions
 \begin{eqnarray}
  \label{l3}
\left.
\begin{array}{c}
 \mathfrak{W}_{+} - \mathfrak{W}_{-} =  0 ,
\\ \partial_y \mathfrak{W}_{+} - \partial_y  \mathfrak{W}_{-} =0
\end{array}
\right\}
\quad \mathrm{when} \ (t,x,y) \in  (0,T) \times \Omega \times \{ 0 \},
\end{eqnarray} 
 then there exist $\lambda_m >0$ and for all $k \in \N$, $\mu_{k,m} > 0$, such that for all $\lambda \geqslant \lambda_m$, 
 \begin{eqnarray}
 \label{esti1}
 | W |_{m,\lambda,T}   \leqslant \frac{\lambda_m}{\lambda}  | f |_{m,\lambda,T}  
\end{eqnarray} 
and for all $ \mu \geqslant \mu_{k,m} $,
 \begin{eqnarray}
   \label{esti2}
    |y^k \, W |_{m,\lambda,T} \leqslant \frac{\mu_{k,m}}{\mu}  \sum_{j=0}^k | y^j \, f |_{m,\mu,T} .
\end{eqnarray} 
\end{proposition}

\begin{proof}

We  multiply the equation (\ref{l1}) by $\mathbf{W}_{\pm}$ and integrate for $(x,y) \in   \Omega \times \R_\pm$.
Hence
\begin{eqnarray}
\label{estim}
(1/2) \partial_t \int_{\Omega  \times \R_\pm} | \mathbf{W}_{\pm} |^2 - J_{1,\pm}- J_{2,\pm}
 =  \int_{\Omega  \times \R_\pm} f_{\pm} . \mathbf{W}_{\pm}
\\ \nonumber \mathrm{where} \ J_{1,\pm} := \int_{\Omega  \times \R_\pm} \mathbf{W}_{\pm} . \partial_y^2 \mathbf{W}_{\pm}
\ \mathrm{and} \ J_{2,\pm} := \int_{\Omega  \times \R_\pm} \mathbf{W}_{\pm} .  (
\mathfrak{W}_{\pm} \wedge \partial_y^2)  \mathbf{W}_{\pm} .
\end{eqnarray} 
 Integrating by parts, we get 
\begin{eqnarray*}
J_{1,\pm} = 
- \int_{\Omega \times \R_\pm } | \partial_y \mathbf{W}_{\pm} |^2
- I_{1,\pm} ,
\ \mathrm{and} \ J_{2,\pm} = 
- \int_{\Omega \times \R_\pm } \mathbf{W}_{\pm} .  ( \partial_y \mathfrak{W}_{\pm} \wedge \partial_y)  \mathbf{W}_{\pm}
- I_{2,\pm} ,
\\ \nonumber 
 \mathrm{where} \ I_{1,\pm} := \int_{\Omega } ( \mathbf{W}_{\pm}. \partial_y \mathbf{W}_{\pm} )|_{y=0} ,
\ \mathrm{and} \   I_{2,\pm} := \int_{\Omega } (\mathbf{W}_{\pm} .  (  \mathfrak{W}_{\pm} \wedge \partial_y  \mathbf{W}_{\pm}))|_{y=0} .
\end{eqnarray*} 
Using the boundary conditions (\ref{l2}) and (\ref{l3}), we get $I_{1,+} - I_{1,-} =  I_{2,+} - I_{2,-} =  0$.
Taking that into account we add the two estimates in (\ref{estim}).
Then we multiply by $e^{-2 \lambda t}$ and integrate in time.
By a Gronwall lemma we get  that there exists $c >0$ such that for all $\lambda \geqslant c$, 
\begin{eqnarray}
\label{estiL2}
 | \partial_y W |^2_{0,\lambda,T} + \lambda | W |^2_{0,\lambda,T} \leqslant c  |<f,W>_{\lambda,T}| .
\end{eqnarray} 

We go on with  estimates tangential to $\{y=0\}$. 
To do this we apply the derivative $\partial_{t,x}^{\alpha} \, $ to the equations (\ref{l1})-(\ref{l2}). 
So we get that $\partial_{t,x}^{\alpha} \, \mathbf{W}_{\pm} $ verify
\begin{eqnarray}
\label{l1tan} 
L ( \mathfrak{W}_{\pm} , \partial_{t},\partial_{y}^2 )
\partial_{t,x}^{\alpha} \, \mathbf{W}_{\pm} = \tilde{f}_{\pm}  \quad \mathrm{when} \ (t,x,y) \in  (0,T) \times \Omega \times \R_\pm ,
\\ \label{l2tan} 
\left.
\begin{array}{c}
\partial_{t,x}^{\alpha} \,  \mathbf{W}_{+} - \partial_{t,x}^{\alpha} \, \mathbf{W}_{-} =  0 ,
\\ \partial_y \partial_{t,x}^{\alpha} \, \mathbf{W}_{+} - \partial_y  \partial_{t,x}^{\alpha} \, \mathbf{W}_{-} =0
\end{array}
\right\}
\quad \mathrm{when} \ (t,x,y) \in  (0,T) \times \Omega \times \{ 0 \},
\end{eqnarray} 
where 
\begin{eqnarray}
\label{source}
\tilde{f}_{\pm} := \partial_{t,x}^{\alpha} \, f_{\pm} 
+ \sum_{|\alpha_1 | + |\alpha_2 | = |\alpha|, |\alpha_2 | < |\alpha|} 
 \partial_{t,x}^{\alpha_1} \, \mathfrak{W}_{\pm}  \wedge  \partial_y^2 \, \partial_{t,x}^{\alpha_2} \, \mathbf{W}_{\pm} .
\end{eqnarray} 
We apply the tangential derivative $\partial_{t,x}^{\alpha} \, $ to the boundary conditions (\ref{l3}) and get
 \begin{eqnarray}
  \label{l3tan}
\left.
\begin{array}{c}
\partial_{t,x}^{\alpha} \, \mathfrak{W}_{+} - \partial_{t,x}^{\alpha} \,  \mathfrak{W}_{-} =  0 ,
\\ \partial_y \partial_{t,x}^{\alpha} \,  \mathfrak{W}_{+} - \partial_y \partial_{t,x}^{\alpha} \,  \mathfrak{W}_{-} =0
\end{array}
\right\}
\quad \mathrm{when} \ (t,x,y) \in  (0,T) \times \Omega \times \{ 0 \},
\end{eqnarray} 
Using the estimate (\ref{estiL2}), we get, for all $\lambda \geqslant c$, 
\begin{eqnarray*}
| \partial_y \partial_{t,x}^{\alpha} \, W |^2_{0,\lambda,T} 
 + \lambda | \partial_{t,x}^{\alpha} \, W |^2_{0,\lambda,T} \leqslant 
 c  |<\tilde{f} , \partial_{t,x}^{\alpha} \, W  >_{\lambda,T}|.
\end{eqnarray*} 
Thanks to (\ref{source}), we get
\begin{eqnarray}
<\tilde{f} , \partial_{t,x}^{\alpha} \, W  >_{\lambda,T} = 
< \partial_{t,x}^{\alpha} \, f  , \partial_{t,x}^{\alpha} \, W  >_{\lambda,T}
+ \sum_{|\alpha_1 | + |\alpha_2 | = |\alpha|, |\alpha_2 | < |\alpha|}  I_{\alpha_1 , \alpha_2 } ,
\end{eqnarray} 
where $ I_{\alpha_1 , \alpha_2 } := I_{+,\alpha_1 , \alpha_2 } + I_{-,\alpha_1 , \alpha_2 }$ with
\begin{eqnarray*}
I_{ \pm,\alpha_1 , \alpha_2 } := < \partial_{t,x}^{\alpha_1} \, \mathfrak{W}_{\pm}  \wedge  \partial_y^2 \,
\partial_{t,x}^{\alpha_2} \, \mathbf{W}_{\pm} , \partial_{t,x}^{\alpha_2} \, \mathbf{W}_{\pm} >_{\lambda,T} .
\end{eqnarray*} 
Using Cauchy-Schwarz inequality, we get
\begin{eqnarray*}
|< \partial_{t,x}^{\alpha} \, f  , \partial_{t,x}^{\alpha} \, W  >_{\lambda,T} | 
\leqslant | f|_{0,\lambda,T} \, .| W|_{0,\lambda,T} .
\end{eqnarray*} 
  We are going to estimate, for all $ \alpha_1 , \alpha_2 $ 
such that $|\alpha_1 | + |\alpha_2 | = |\alpha|, |\alpha_2 | < |\alpha|$, the term $ I_{\alpha_1 , \alpha_2 }$.
Integrating by parts, we get
$I_{ \pm,\alpha_1 , \alpha_2 } :=  \sum_{l=1}^3 \, I^l_{ \pm,\alpha_1 ,
\alpha_2 },$
with
\begin{eqnarray*}
 I^1_{ \pm,\alpha_1 , \alpha_2 } &:=& - 
 < \partial_{t,x}^{\alpha_1} \, \partial_y \mathfrak{W}_{\pm}  \wedge  \partial_y \, \partial_{t,x}^{\alpha_2} \, \mathbf{W}_{\pm}
 , \partial_{t,x}^{\alpha} \, \mathbf{W}_{\pm} >_{\lambda,T} ,
\\  I^2_{ \pm,\alpha_1 , \alpha_2 } &:=& -  
 < \partial_{t,x}^{\alpha_1} \,  \mathfrak{W}_{\pm}  \wedge  \partial_y \, \partial_{t,x}^{\alpha_2} \, \mathbf{W}_{\pm} ,
\partial_{t,x}^{\alpha} \, \partial_y \mathbf{W}_{\pm}>_{\lambda,T} ,
\\ I^3_{ \pm,\alpha_1 , \alpha_2 } &:=& \mp   <<
\{  (\partial_{t,x}^{\alpha_1} \, \mathfrak{W}_{\pm} \wedge  \partial_y \, \partial_{t,x}^{\alpha_2} \, \mathbf{W}_{\pm} )\}|_{y=0},
\{ \partial_{t,x}^{\alpha} \, \partial_y \mathbf{W}_{\pm} \}|_{y=0} >>_{\lambda,T} ,
\end{eqnarray*} 
where $<<.,.>>_{\lambda,T}$ denotes the scalar product of $L^2 ((0,T) \times \Omega )$ associated to the mesure $e^{-\lambda t} dt
dx$. 
 Thanks to the boundary conditions   (\ref{l2tan})-(\ref{l3tan}), 
we get $ I^3_{ +,\alpha_1 , \alpha_2 } - I^3_{ -,\alpha_1 , \alpha_2 } = 0$.
Using Cauchy-Schwarz inequality, we get
\begin{eqnarray*}
| I^1_{ \pm,\alpha_1 , \alpha_2 }| &\leqslant& 
| \partial_{t,x}^{\alpha_1} \, \partial_y \mathfrak{W}_{\pm}  \wedge  \partial_y \, \partial_{t,x}^{\alpha_2} \, \mathbf{W}_{\pm}
|_{0,\lambda,T} \, . ||  \mathbf{W}_{\pm}||_{m,\lambda,T} ,
\\ | I^2_{ \pm,\alpha_1 , \alpha_2 }| &\leqslant& |  \partial_{t,x}^{\alpha_1} \,  \mathfrak{W}_{\pm}  \wedge  \partial_y \,
\partial_{t,x}^{\alpha_2} \, \mathbf{W}_{\pm} |_{0,\lambda,T} \, . || \partial_y  \mathbf{W}_{\pm}||_{m,\lambda,T} ,
 \end{eqnarray*} 
Using Gargliardo-Nirenberg inequalities, we get
\begin{eqnarray*}
| I^1_{ \pm,\alpha_1 , \alpha_2 } |   &\leqslant 
 & \quad c ( ||\partial_y \mathfrak{W}_{\pm} ||_{m,\lambda,T} . ||  \mathbf{W}_{\pm}   ||_{Lip} 
 +   || \mathfrak{W}_{\pm} ||_{Lip} . || \partial_y \mathbf{W}_{\pm} ||_{m,\lambda,T} ). ||\mathbf{W}_{\pm}||_{m,\lambda,T} ,
\\  | I^2_{ \pm,\alpha_1 , \alpha_2 } | &\leqslant 
 & \quad c ( || \mathfrak{W}_{\pm} ||_{m-1,\lambda,T} . ||  \mathbf{W}_{\pm}   ||_{Lip} 
 +   | \mathfrak{W}_{\pm} |_{Lip} . || \partial_y \mathbf{W}_{\pm} ||_{m-1,\lambda,T} ). || \partial_y  \mathbf{W}_{\pm}||_{m,\lambda,T} .
\end{eqnarray*} 
Hence we get 
\begin{eqnarray*}
|I_{\alpha_1 , \alpha_2 } | &\leqslant& \frac{1}{2} || \partial_y  \mathbf{W}_{\pm}||^2_{m,\lambda,T} 
+ C ( || \mathbf{W}_{\pm}||^2_{m,\lambda,T} + || \partial_y  \mathbf{W}_{\pm}||^2_{m-1,\lambda,T}).
\end{eqnarray*} 
We deduce that there exists $\lambda_m >0$ 
such that for all $\lambda \geqslant \lambda_m$, there holds 
$ || W ||_{m,\lambda,T}  \leqslant \frac{\lambda_m}{\lambda}  || f
||_{m,\lambda,T} $.

To prove  the estimates (\ref{esti1}), it remains to get 
 normal estimates.
The cases $\alpha_4 =0$ or $1$ are already treated in the tangential estimates.
If  $\alpha_4 \geqslant 2$, we proceed by iteration, extirpating $\partial_y^2 \mathbf{W}_{\pm}$ from the equations.

It remains to get the estimates (\ref{esti2}). 
First we notice that  for  $p \geqslant 1$ the function $y^p \, \mathbf{W}_{\pm}$ verify the initial boundary value problem
\begin{eqnarray*}
L (\mathfrak{W}_{\pm}, \partial_{t}, \partial_{y}^2 )
\mathbf{W}_{\pm}^{[p]} = f_{\pm}^{[p]} \quad \mathrm{when} \ (t,x,y) \in  (0,T) \times \Omega \times \R_\pm ,
\\ \left.
\begin{array}{c}
 \mathbf{W}^{[p]}_{+} - \mathbf{W}^{[p]}_{-} =  0 ,
\\ \partial_y \mathbf{W}^{[p]}_{+} - \partial_y  \mathbf{W}^{[p]}_{-} =0
\end{array}
\right\}
\quad \mathrm{when} \ (t,x,y) \in  (0,T) \times \Omega \times \{ 0 \},
\\  \mathbf{W}^{[p]}_{\pm} = 0 \quad \mathrm{when} \ (t,x,y) \in  \{ 0 \} \times \Omega \times \R_\pm  ,
\end{eqnarray*} 
where
\begin{eqnarray*}
f_{\pm}^{[p]} &=& y^p \, f_{\pm} 
+ \sum_{j=0}^{p-1} \, ( q^1_j \partial_y  \mathbf{W}^{[j]}_{\pm} +  q^2_j \, \mathfrak{W}_{\pm}  \wedge \partial_y 
\mathbf{W}^{[j]}_{\pm}),
\end{eqnarray*} 
where the $q^1_j$ and the $q^2_j$ are in $\N$. 
Thus we prove, by iteration on $p$  and thanks to the inequality
(\ref{esti1}),  the estimate
\begin{eqnarray*}
\sqrt{\mu} || \partial_y (y^p \, W)||_{m,\mu,T} + \mu || y^p \, W||_{m,\mu,T}  \leqslant \sum_{j=0}^p || y^j \,
f||_{m,\mu,T} 
\end{eqnarray*} 
which implies the estimate (\ref{esti2}). 

\end{proof}
\begin{stepbd}
We use an iterative scheme.
\end{stepbd}
We define the iterative scheme $(\mathbf{W}^\nu_{\pm})_{\nu \in \N}$ by setting
$ \mathbf{W}^0_{\pm}$ equal to zero 
and, by iteration, when $\mathbf{W}^\nu_{\pm}$ is defined, we take $\mathbf{W}^{\nu+1}_{\pm}$ as solution of 
\begin{eqnarray*}
  L (V_{\pm} + \mathbf{W}^\nu_{\pm} , \partial_{t}, \partial_{y}^2 )
\mathbf{W}^{\nu+1}_{\pm} = \hat{F}(t,x,y,  \mathbf{W}^\nu_{\pm} , \partial_y
\mathbf{W}^\nu_{\pm}) \quad  \mathrm{when} \ (t,x,y) \in  (0,\infty) \times \Omega \times \R_\pm ,
\\  \left.
\begin{array}{c}
 \mathbf{W}^{\nu+1}_{+} - \mathbf{W}^{\nu+1}_{-} =  0 ,
\\ \partial_y \mathbf{W}^{\nu+1}_{+} - \partial_y  \mathbf{W}^{\nu+1}_{-} =0
\end{array}
\right\}
\quad \mathrm{when} \ (t,x,y) \in  (0,T) \times \Omega \times \{ 0 \} ,
\\ \mathbf{W}^{\nu+1}_{\pm} = 0 \quad \mathrm{when} \ (t,x,y) \in  \{ 0 \} \times \Omega \times \R_\pm .
\end{eqnarray*} 

Thanks to the linear estimates, to a Sobolev embedding and to some Gargliardo-Nirenberg inequalities, we show that the
iterative scheme $(\mathbf{W}^\nu_{\pm})_{\nu \in \N}$ converge, when $\nu \rightarrow +\infty $  toward 
 some  solutions $\mathbf{W}_{\pm}  \in \mathcal{N}_{\pm} (T)$ of the
 problem (\ref{++})-(\ref{np2}). By going back to the original problem
 (\ref{p1})-(\ref{p2}), the first sentence of Theorem $\ref{profil0}$ is
 now proved. When $x \notin \mathcal{V}_\Sigma$, the function $u^0_+ -
 u^0_-$ in the right hand side of  (\ref{p1}) vanishes and so do the
 functions  $\mathcal{U}_{\pm}  $.

\end{proof}
\begin{remark}
\rm
Notice that the possibility of a blow-up  can be controlled with Lipschitz norm in a very classical way. 
However we do not know whether the solutions $\mathcal{U}$ actually
blow-up or globally exist.
\end{remark}
\subsection{Construction of  $\mathfrak{U}$}
\label{deja}

In this section we define the  boundary layer profile
$\mathfrak{U}$ as a solution of a linear boundary value problem. 
Let us recall that this function describes a boundary layer which appears near the boundary 
to compensate the lost of the Neumann condition
 from the complete model (\ref{LL1})-(\ref{LL2})-(\ref{LL3}) to
the limit model  (\ref{LL0})  ($\eps = 0$). 
Such a boundary layer was already mentioned in paper \cite{CFG1}.
Let $\Theta$  be a $C^\infty$  function on
$\Omega$ such that $\Theta =1$   in a neighborhood $\mathcal{W}_\Gamma$
 of $\Gamma$ such that
$\mathcal{W}_\Gamma \subset \subset  \mathcal{V}_\Gamma$
 and  $\Theta =0$  in $\Omega -
 \mathcal{W}_\Gamma$.

\begin{theorem}
There exists $\mathfrak{U}  \in \mathcal{N}_{+} (T)$ which verifies
\begin{eqnarray*}
 L ( u^0 , \partial_{t}, \partial_{z}^2 ) \mathfrak{U}
& = & -( \mathfrak{U}. n) u^0 \wedge   n
 + \mathfrak{U}  \wedge   \cH(u^0 ) 
\\ \nonumber && + \mathfrak{U} \wedge ( u^0 \wedge \cH(u^0 ) )
-( \mathfrak{U}. n)   u^0 \wedge ( u^0  \wedge  n )
 +u^0 \wedge ( \mathfrak{U}  \wedge \cH(u^0 ) ,
 \end{eqnarray*}
when $(t,x,z) \in  (0,T) \times \Omega \times \R_+ $,
\begin{eqnarray}
  \partial_z \mathfrak{U}  = \Theta (x) \partial_\mathfrak{n} u^0
  \quad \mathrm{when} \ (t,x,z) \in  (0,T) \times \Omega \times \{ 0 \}.
\end{eqnarray}
Moreover there holds $\mathfrak{U}(t,x,z)=0$ for $x  \notin \mathcal{V}_\Sigma$.
\end{theorem}

\begin{proof}
Proceeding as in the proof of Theorem $\ref{profil0}$, we prove the existence of compatible initial data. 
Then we follow the proof of Proposition $4.2$ of \cite{CFG1}.
\end{proof}

\subsection{Construction of $\mathbf{w}^\eps$ }
\label{s3}
In this section, we look at the remainder $\mathbf{w}^\eps$.
We will proceed in four steps.
First in section \ref{c1} we  will reduce the initial problem  (\ref{LL1})-(\ref{LL2})-(\ref{LL3})
for the unknown  $u^\eps$ to a problem for $\mathbf{w}^\eps$. 
Indeed in order  to get a homogeneous boundary problem,  we will add a
corrector to $\mathbf{w}^\eps$ and rather work with the resulting term $w^\eps$. 
Some Borel classical arguments will insure the existence of convenient
initial data for the resulting reduced problem which means that
compatibility conditions either on $\Gamma$ and on $\Sigma$ are satisfied.
We will  prove that the solutions of this nonlinear problems
exist not only for a common non trivial time, in fact even till the
lifetime $T$ of the profiles $\mathcal{U}$.
Moreover these  solutions satisfy some estimates uniform with respect to
$\eps$.
 The method lies on  a simple Picard iterative scheme  (cf. section
 \ref{c2}) and on linear estimates (cf. section
 \ref{c3}). 
More precisely we will use  $L^2$-type conormal estimates of  only the
two first normal derivatives, and some Lipschitz estimates.
A few carefulness reveals that the presence of the operator $\mathcal{H}$ does not
cause any loss of factor $\eps$ or any loss of derivatives.

\subsubsection{A reduced problem}
\label{c1}

Since we look for  solutions  $u^\eps$ of (\ref{LL1})-(\ref{LL2})-(\ref{LL3})
 of the form (\ref{decomp}) where the functions 
 \begin{eqnarray*}
\label{decomp2}
a^{\eps}  (t,x) &:=&  \mathcal{U}  (t,x, \frac{\Psi(x)}{\eps}) 
+  \eps \Big(\mathfrak{U} (t,x, \frac{\Phi(x)}{\eps} )  +  \mathbf{w}^{\eps}  (t,x)  \Big)
\end{eqnarray*} 
 have been constructed above, 
we look for a problem in term of the  remainder $\mathbf{w}^\eps$.
In fact, in order to get a homogeneous boundary problem,  we choose a function $\rho(t,x) \in H^\infty$ such that
\begin{eqnarray}
\label{B1}
  \Dn \rho\, |_\Gamma  &=& - \Dn   \mathfrak{U} (t,x,0)|_{\Gamma}.
 \end{eqnarray}
 and will look for remainders $\mathbf{w}^\eps$
 of the form $\mathbf{w} ^\eps =  \rho + w^\eps $. 
 Let us explain why. 
On the boundary $\Gamma$, the function  $a^\eps $ satisfies:
\begin{equation}
\label{B2}
\Dn a^\eps|_{\Gamma} = \eps \, \Dn \mathfrak{U} (t,x,0)|_{\Gamma},
\end{equation}
Hence in general $a^\eps$ does not satisfy the homogeneous Neumann
boundary condition on $  \Gamma$. 
 We  define the function $ \ta^\eps := a^\eps + \eps \rho$.
Thus we look for solutions  $u^\eps$ of (\ref{LL1})-(\ref{LL2})-(\ref{LL3})
 of the form $ u^\eps =   a^\eps + \eps \, \mathbf{w}^\eps     =    \ta^\eps + \eps \, w^\eps $. 
 Combine (\ref{LL2}), (\ref{B1}) and (\ref{B2}) to find a homogeneous
 Neumann boundary condition on $\Gamma $ for $w^\eps$:
\begin{eqnarray}
\label{5}
\Dn w^\eps  &=& 0  \quad  \mathrm{on} \ ]0,T[ \times \Gamma.
\end{eqnarray}

We now look for an equation on the unknown $w^\eps$.
The function $\ta^\eps$ belongs to  $\cC^1((0,T)\times \Omega)$ and
to $H^\infty_\Sigma(\Omega)$.
Moreover, $\ta^\eps$ satisfies the equation
\begin{equation}\label{2}
    \cL(\ta^\eps,\D)\, \ta^\eps =
    \bF\big(\ta^\eps, \eps \D_x \ta^\eps,
     \cH(\ta^\eps)\big) + \eps  r^\eps
\end{equation}
where the family $(r^\eps )_\eps$ lies in the set $E$ (defined above
Theorem $\ref{main2}$).
The system for the unknown $w^\eps(t,x)$ writes
\begin{equation}
    \begin{aligned}
    \label{4} \cL ( \ta^\eps + \eps w^\eps  , \D) w^\eps =
K(\eps,\ta^\eps,&\eps \D_x \ta^\eps,\cH(\ta^\eps),
w^\eps,\eps\D_x w^\eps,\cH(w^\eps)) +  r^\eps
 \quad  \mathrm{in } \  ]0,T[ \times \Omega
    \end{aligned}
\end{equation}
where $K$ is a smooth function of its arguments. Let us use more
concise notations, and note
\begin{eqnarray}
A^\eps := \big(\, \ta^\eps,\eps \D_x
\ta^\eps,\cH(\ta^\eps) \, \big)
\quad \text{and} \quad 
W^\eps := \big(\,  w^\eps,\eps\D_x w^\eps,\cH(w^\eps)  \, \big).
\end{eqnarray}
Then, the Taylor formula shows that the function $K$ has the
following form:
\begin{equation}
\nonumber
K(\eps,A^\eps, W^\eps) = G(\eps, A^\eps, \eps  W^\eps) W^\eps
\end{equation}
where $G$ depends smoothly on its arguments (including $\eps$),
which will be useful in the sequel.

 Following \cite{S1} there exist a family $(w^\eps_{\text{init}})_\eps$ 
 of  compatible  initial conditions for the problem (\ref{4})-(\ref{5})
 which verifies suitable uniform estimates with respect to $\eps$.
We choose such a family.

\subsubsection{The iterative scheme}
\label{c2}

We want to solve the problem (\ref{4}),(\ref{5}). We use
a simple Picard(-Banach-Caccioppoli)  iterative scheme defining a sequence $w^{\eps,\nu}$
which will converge to the solution of the problem. For clarity, we
adopt the following more concise notations
$$
A^\eps :=\big(\, \ta^\eps,\eps \D_x
\ta^\eps,\cH(\ta^\eps)\, \big) \quad    \text{and}   \quad   W^{\eps,\nu} :=
\big(\, w^{\eps,\nu},\eps\D_x w^{\eps,\nu}, \cH(w^{\eps,\nu}) \,
\big).
$$
With these notations, the iterative scheme writes
\begin{equation}
    \begin{aligned}
    \label{50} \cL ( \ta^\eps + \eps w^{\eps, \nu}  , \D)
    w^{\eps,\nu+1} = f^{\eps, \nu}
 \quad  \mathrm{in } \  ]0,T[ \times \Omega
    \end{aligned}
\end{equation}
where
\begin{equation}\label{91}
   f^{\eps,\nu} := G(\eps,A^\eps,
\eps W^{\eps,\nu}) W^{\eps,\nu} +  r^\eps
\end{equation}
This equation
is coupled with the initial and boundary conditions:
\begin{eqnarray}\label{51}
  \Dn w^{\eps, \nu +1} &=& 0 \ \mathrm{on} \ ]0,T[ \times \Gamma \\
  w^{\eps, \nu +1}|_{t=0} &=& w^\eps_{\text{init}}.\label{52}
\end{eqnarray}
The iterative scheme  is initialized with $w^{\eps,0}(t,x) :=
w^\eps_{\text{init}}(x)$.

\subsubsection{Estimates for a linear parabolic system}
\label{c3}

Consider the linear problem
\begin{eqnarray}
\label{num1}
  \cL(\ta^\eps + \eps \bb, \D) \bu  &=& f \
  \mathrm{on} \ ]0,T[ \times \Omega\\
\label{num2}  \Dn \bu &=& 0 \ \mathrm{on} \ ]0,T[ \times \Gamma ,
\end{eqnarray}
We endow the space $H^m_{co}(]0,T[\times \Omega)$  with the usual weighted norm with $\lambda \geq 1$:
$$
\|  \bu \|_{m,\lambda} := \sum_{|\alpha| \leq m \, , \, \alpha \in
\NN^{1+\mu}} \lambda ^{m-|\alpha|}\| e^{-\lambda t} \cZ^\alpha \bu
\|_{L^2(]0,T[ \times \Omega)}.
$$
 In order to estimate the initial data, we introduce the similar norms
built with the set $\cT_0$ instead of $\cT$, integrating on $\Omega$
instead of $[0,T]\times \Omega$:
$$
| \bu |_{m,\lambda} := \sum_{|\alpha| \leq m \, , \, \alpha_0 = 0\, ,
\, \alpha \in \NN^{1+\mu}} \lambda ^{m-|\alpha|}\|  \cZ^\alpha \bu
\|_{L^2(\Omega)}.
$$
We will use the following classical
Gagliardo-Moser-Nirenberg estimates  for conormal derivatives (see
\cite{G2}).
\begin{lemma}\label{moser}
Let $m\in \NN$. There is $c_m>0$ such that, for any $a_1, \dots, a_k
\in H^m_{co}(]0,T[\times \Omega) \cap L^\infty(]0,T[\times \Omega)$,
 for all multi-index $\alpha_1 \in \NN^{\mu +1}, \dots , \alpha_k
\in \NN^{\mu +1}$, with $|\alpha_1|+\dots +|\alpha_k| \leq m$, for all $\lambda \geq 1$:
\begin{equation}\label{89}
\| \cZ^{\alpha_1} a_1 \dots \cZ^{\alpha_k}a_k \|_{0,\lambda} \ \leq
\ c_m \ \sum_{1 \leq j\leq k} \Big( \| a_j\|_{m,\lambda} \prod_{i
\ne j} \|a_i\|_\infty \Big).
\end{equation}
\end{lemma}

The following proposition gives some $\eps$-conormal estimates for the
two first normal derivatives of the solutions of the problem 
(\ref{num1})-(\ref{num2}).

\begin{proposition}
\label{22}
Let $R>0$ be an arbitrary constant and $m \geqslant 3$.  There exist
$C_m(R)>0$ and $\lambda_m>0$ such that for $\sigma$ fixed constant large enough, depending only on
the choices of the vector fields $\cZ_j$, the following holds true. Assume that
\begin{equation}\label{44}
\eps  \, ( \ \|\bb\|_\infty +\sum_{0\leq j\leq \mu} \|\cZ_j
\bb\|_\infty+\|\eps \D_x \bb\|_\infty \ ) \, \leq \, R,
\end{equation}
then, for all $\lambda \geq \lambda_m$,  the following estimates hold:
\begin{equation}\label{20}
\begin{aligned}
   \|\eps  \D_x \bu\|_{m,\lambda} +
   \lambda \|\bu\|_{m,\lambda}
    \ &\leq \ C_m(R) \ \big[ \ \lambda^{-1}\ \| f \|_{m,\lambda} + \ I_{m,\lambda}(\bu)\\
    &+ \ \eps \, (\, \|\eps \D_x \bb \|_{m,\lambda} +
    \| \bb\|_{m,\lambda}\, )\ (\, \|\bu\|_\infty + \|\eps\D_x \bu\|_\infty \, )\ \big],
    \end{aligned}
\end{equation}
where
$$
I_{m,\lambda}(\bu) := \sum _{0 \leq k \leq m}|(\D_t^k \bu)_{|t=0}
|_{m-k,\lambda}.
$$
and
\begin{equation}\label{43}
\begin{aligned}
\|(\eps \Dn)^2 \bu \|_{m,\lambda} \ \leq \ C_{m}(R) \big[ \
\|f\|_{m,\lambda} + \| \bu \|_{m+1,\lambda} 
+ \eps \|\bb\|_{m +1,\lambda}\,(\|\bu\|_{\infty} +
\|f\|_\infty) \  \ \\
+ \eps \| \eps \Dn \bu \|_{m +1,\lambda} + \eps^2 \| \bu
\|_{m+2,\lambda} \ \big].
\end{aligned}
\end{equation}
\end{proposition}

\begin{proof}

\emph{Step 1}.
Let us note $\bv:= e^{-\lambda t} \, \bu$, which satisfies
\begin{eqnarray}
  \cL(\ta^\eps_{app} + \eps \bb, \D) \bv + \lambda \bv &=&
   e^{-\lambda t} f \label{10}\
  \mathrm{on} \ ]0,T[ \times \Omega\\
  \Dn \bv &=& 0 \ \mathrm{on} \ ]0,T[ \times \Gamma.\\
  \bv &=& w^\eps_{\text{init}} \ \mathrm{on} \ t=0.
\end{eqnarray}
 Let us note $\| . \|_{L^2}$ the $L^2$ norm in
 $[0,T]\times \Omega$, and $|.|_{L^2}$
 the $L^2$ norm in $\Omega$. Multiplying  (\ref{10}) by $\bv$ and integrating 
on $]0,T[ \times \Omega$ gives the following estimate,
integrating by parts
the $\eps^2 \Delta_x$ with Green's formula in $\Omega$:
\begin{equation}\label{11}
    \eps^2 \|\nabla_x \bv \|^2_{L^2}+ \lambda \| \bv \|^2_{L^2}
    \leq  2 \ |(( e^{-\lambda t}f, \bv ))_{L^2}|
    + |\bv(0)|_{L^2} ,
\end{equation}
for all $\lambda \geq \lambda_0$ if $\lambda_0$ is fixed
large enough, and for all $\eps >0$. In terms of $\bu$ it writes
\begin{equation}\label{11'}
    \eps^2 \|\nabla_x \bu \|^2_{0,\lambda}+ \lambda \| \bu \|^2_{0,\lambda}
    \leq  2 \ |(( f, \bu ))_{L^2_\lambda}|+ |\bu(0)|_{L^2},
\end{equation}
where $L^2_\lambda $ is the Hilbert space $L^2(]0,T[ \times \Omega, d\mu)$
with the measure $d\mu := e^{-2\lambda t} dt dx$.

Using now the Cauchy-Schwarz inequality in the right hand side, and absorbing in the left hand side the term
in $\|v\|_{L^2}^2$ yields the desired estimate
for $m=0$ and some constant $c_0 >0$.
\\
\\ \emph{Step 2}. We show the inequality by induction on $m$.
Assume it for $m-1$. We apply a tangential operator $\cZ^\alpha$
with fields $\cZ_i \in \cT$ to the system, and $|\alpha| = m$. The function
$\cZ^\alpha \bu$ satisfies the same boundary conditions.
The $L^2$ estimate (\ref{11'}) gives, for $\lambda \geq \lambda_0$:
\begin{equation}\label{21}
\begin{aligned}
     \eps^2 \|\nabla_x \cZ^\alpha \bu \|^2_{L^2}+ \lambda \| \cZ^\alpha \bu \|^2_{L^2}
    \leq&  2 |(( e^{-\lambda t}\cZ^\alpha f  +  [(\ta^\eps_{app} + \eps \bb)\eps^2\Delta_x,\cZ^\alpha ]\wedge
    \bu   , \cZ^\alpha\bu ))_{L^2_\lambda}|.
\end{aligned}
\end{equation}
where $[.,.]$  denotes the commutator. 
Using Cauchy-Schwarz inequality and $2ab\leq 2\lambda^{-1} \, a^2 + \lambda b^2/2$ yields:
\begin{equation}\label{12}
\begin{aligned}
     \eps^2 \|\nabla_x \cZ^\alpha \bu \|^2_{L^2}+ \frac{\lambda}{2} \| \cZ^\alpha \bu \|^2_{L^2}\
    \leq&  \ \frac{2}{\lambda}\ \|e^{-\lambda t}\cZ^\alpha f\|^2_{L^2}\\
    + \ 2 \ |(([(\ta^\eps_{app} &+ \eps \bb)\eps^2\Delta_x,\cZ^\alpha ]\wedge \bu , \cZ^\alpha \bu))_{L^2_\lambda}|.
\end{aligned}
\end{equation}
We need to control the second term in the right hand side of (\ref{12}). The commutator
$[\ \ta^\eps_{app}\, \eps^2\Delta_x,\cZ^\alpha ]$
writes as a finite sum
\begin{equation}\label{13}
    \eps ^2  \sum_{|\beta|\leq m+1} a_\beta^\eps(t,x) \cZ^\beta +
    \eps   \sum_{|\gamma|\leq m} b_\gamma^\eps(t,x) \eps \Dn \cZ^\gamma
    +   \sum_{|\delta|\leq m-1} c_\delta^\eps(t,x) (\eps\Dn)^2
     \cZ^\delta
\end{equation}
where the coefficients $a_\beta^\eps$ , $b_\gamma^\eps$,
$c_\delta^\eps$ are bounded functions satisfying
\begin{equation}
\label{14}    \sup _{\eps \in ]0,1]}\|\eps \Dn a_\beta^\eps\|_{L^\infty(\Omega)} 
    + \|\eps \Dn b_\gamma^\eps\|_{L^\infty(\Omega)} + \|\eps \Dn c_\delta^\eps \|_{L^\infty(\Omega)} < \infty
\end{equation}
for all $\beta, \gamma, \delta$, because (\ref{14})
holds clearly if we replace  $L^\infty(\Omega)$ by
 $L^\infty(\Omega_+)$ or by $L^\infty(\Omega_-)$,
 and because $\ta^\eps_{app}$ is in
$H^1(\Omega)$ for all $\eps >0$.
Hence we are led to control the corresponding three sort of terms:
\begin{eqnarray}
\label{15}
\eps^2 (( \ a_\beta^\eps \cZ^\beta \bu \, ,
    \, \cZ^\alpha \bu \ ))_{L^2_\lambda} ,  \
    \eps (( \ b_\gamma^\eps (\eps \Dn) \cZ^\gamma \bu \, ,
    \, \cZ^\alpha \bu \ ))_{L^2_\lambda} , \
    (( \ c_\delta^\eps (\eps \Dn)^2 \cZ^\delta \bu \, ,
    \, \cZ^\alpha \bu \ ))_{L^2_\lambda},
\end{eqnarray}
where $|\beta| \leq m+1$, $|\gamma| \leq m$, $|\delta| \leq m-1$.
The first two terms in (\ref{15}) are simply controlled by
$\delta  \| \eps \nabla_x \bu \|^2_{m,\lambda} +
C_\delta \, \delta^{-1} \,  \| \bu \|^2_{m,\lambda}$
for $\delta$ arbitrarily small, and $C_\delta$ being a constant
depending on $\delta$, but independent of $\eps$.
For the third term one uses an integration by parts (by Green's formula) of the field $\Dn$ to show that this term writes as a sum of terms of the form
$$
d^\eps \
\eps^{2-j-j'}(( (\eps \Dn)^j\cZ^{\delta}  \bu, (\eps \Dn)^{j'} \cZ^{\alpha} \bu ))_{L^2_\lambda}
$$
where $|\delta| \leq m-1$, $j,j'\in \{0,1\}$,
and $d^\eps$ is a bounded function (uniformly in $\eps$)
since all the boundary terms terms vanishes:
$\Dn \cZ^\alpha \bu|_{\D \Omega} = 0$,  for all $\alpha \in \RR^\mu.$
It follows that the third term in (\ref{15}) is controlled by
$C\lambda^{-1} \| \eps\nabla_x \bu \|^2_{m,\lambda}
+ C\| \bu \|^2_{m,\lambda}$
for a constant $C$ independent of $\eps$, and all $\lambda \geq 1$.
Hence, by choosing  a $\delta >0$ arbitrarily small, and $\lambda_1>0$ large enough, there holds
\begin{equation}
\nonumber
     | \ (( \ [\ta^\eps_{app} \eps^2\Delta_x,\cZ^\alpha ]\wedge \bu \, ,
    \, \cZ^\alpha \bu \ ))_{L^2_\lambda} \ | \ \leq \
    \delta  \| \eps \nabla_x \bu \|^2_{m,\lambda} +
c_m \| \bu \|^2_{m,\lambda}
\end{equation}
for all $\lambda \geq \lambda_1$, and for all $\eps\in ]0,1]$,
with a constant $c_m$ independent of $\eps$.

We need now to estimate the term
\begin{equation}\label{24}
(( [\ \eps \bb\, \eps^2\Delta_x,\cZ^\alpha ]\wedge \bu , \cZ^\alpha \bu))_{L^2_\lambda}.
\end{equation}
The commutator $[\ \bb\, \eps^2\Delta_x,\cZ^\alpha ]$
writes as a finite sum
\begin{eqnarray}
\nonumber
\eps ^2  \sum_{|\beta |\leq m,|\beta '|\leq m+1, |\beta| + |\beta'| \leq m+2} a_{\beta,\beta '}
(\cZ^\beta\bb )\cZ^{\beta'} \label{25}\\
\nonumber
   + \ \eps   \sum_{|\gamma |\leq m,|\gamma '|\leq m, |\gamma|+|\gamma'| \leq m+1} b_{\gamma,\gamma '}(\cZ^{\gamma }\bb )(\eps \Dn) \cZ^{\gamma '} \label{26}\\
 \nonumber
   +  \  \sum_{|\delta |\leq m,|\delta '|\leq m-1, |\delta| + |\delta'| \leq m } c_{\delta, \delta '}
    (\cZ^{\delta }\bb)(\eps\Dn)^2
     \cZ^{\delta '} \label{27}
\end{eqnarray}
where $a_{\beta,\beta '}, b_{\gamma,\gamma '},
c_{\delta, \delta '}$ are smooth fonctions on $\overline{\Omega}$. Hence to control the term (\ref{24}) we are led to
estimate tri-linear terms in $(\bb,\bu,\bu)$ of the following form
(where $d\mu := e^{-2\lambda t}dt dx$):
\begin{eqnarray}
\eps ^2 \int _{]0,T[\times \Omega} a_{\beta, \beta '}
\cZ^\beta\bb .\cZ^{\beta'}\bu_i . \cZ^\alpha \bu_j \ d\mu ,
\quad |\beta|\leq m, |\beta '|\leq m+1, |\beta| + |\beta'| \leq m+2 \label{28}\\
\eps \int _{]0,T[\times \Omega} b_{\gamma,\gamma '}\ \cZ^{\gamma }\bb \ . \ \eps \Dn \cZ^{\gamma '}\bu_i \ . \ \cZ^\alpha \bu_j \ d\mu , \quad |\gamma|\leq m, |\gamma '|\leq m, |\gamma|+|\gamma'| \leq m+1 \label{29}\\
\int _{]0,T[\times \Omega} c_{\delta,\delta '}\ \cZ^{\delta }\bb . \,
(\eps \Dn)^2 \cZ^{\delta '}\bu_i \, .  \cZ^\alpha \bu_j \ d\mu , \quad
|\delta|\leq m, |\delta '|\leq m-1, |\delta|+|\delta'| \leq m  ,   \label{30}
\end{eqnarray}
where the $\bu_i$ are the components of the vector $\bu$.
Let us treat the term (\ref{30}). By the green formula, the integral can be written as a sum of integrals of the form
\begin{eqnarray}
     \int _{]0,T[\times \Omega} c_{\delta,\delta '}\ \cZ^{\delta }\eps \Dn\bb . \cZ^{\delta '}\eps \Dn\bu_i .  \cZ^\alpha \bu_j \ d\mu \label{31}\\
    \int _{]0,T[\times \Omega} c_{\delta,\delta '}\ \cZ^{\delta }\bb . \, \cZ^{\delta'}\eps\Dn\bu_i .\,  \cZ^\alpha\eps\Dn\bu_j \ d\mu ,\label{32}\\
    \eps \ \int _{]0,T[\times \Omega} d_{\delta,\delta '}\ \cZ^{\delta }\bb . \, \cZ^{\delta '}\eps \Dn\bu_i \, . \cZ^\alpha \bu_j \ d\mu ,\label{33}
\end{eqnarray}
and other terms involving lower order derivatives easy to control.
The term (\ref{31}) is controlled by
$$
c \| \eps \Dn \cZ^\delta \bb \eps \Dn \cZ^{\delta '}\bu_i\|_{0,\lambda}
\ \| \bu_j \|_{m,\lambda},
$$
which is bounded by using the Gagliargo-Nirenberg-Moser estimate by
\begin{equation}
\nonumber
    c\big( \ \| \eps \Dn \bb \|_{m,\lambda} |\eps \Dn \bu|_\infty +
    \| \eps \Dn \bu \|_{m,\lambda} |\eps \Dn \bb|_\infty \ \big)\
    \|\bu\|_{m,\lambda}
\end{equation}
and hence by
\begin{equation}
\nonumber
    c(1+R) \big( \ \| \eps \Dn \bb \|_{m,\lambda} |\eps \Dn \bu|_\infty +
    \| \eps \Dn \bu \|_{m,\lambda}  \big)\
    \|\bu\|_{m,\lambda}.
\end{equation}
For the term (\ref{32}) there are two cases. The first case is when $\delta = 0$.
In that case the integral is bounded by
$$
c \ \| \eps \Dn \bu_i \|_{m-1,\lambda} \
\| \bu_j \|_{m, \lambda} \ \leq \
\lambda ^{-1} \| \eps \Dn \bu \|_{m,\lambda} ^2.
$$
The second case is when $|\delta| \geq 1$. In that case we write
$\cZ^\delta \bb = \cZ^{\delta"} \cZ_k\bb$ and apply the Gagliardo-Nirenberg-Moser inequality with $\cZ\bb$ in $L^\infty$.
The term in bounded by
$$
c\big( \ \| \cZ\bb \|_{m-1,\lambda} |\eps \Dn \bu|_\infty +
    \| \eps \Dn \bu \|_{m-1,\lambda} |\cZ \bb|_\infty \ \big)\
    \|\eps \Dn \bu\|_{m,\lambda}
$$
and hence by
\begin{equation}
\nonumber
  c \ \| \bb \|_{m,\lambda} |\eps \Dn \bu|_\infty  \|\eps \Dn \bu\|_{m,\lambda}
    +
    c R \lambda ^{-1}  \| \eps \Dn \bu \|_{m,\lambda}^2
\end{equation}
The next terms like (\ref{33}) are easier to treat in the same way,
and are bounded by the same terms.
The term (\ref{30}) was the more delicate to
estimate. The terms (\ref{29}) and (\ref{28}) are
simpler and can be treated in a similar way.
Replacing in the right hand side of (\ref{21}) and summing over all the possible
operators $\cZ^\alpha$ gives the desired estimate, and the proposition
is proved.
\end{proof}

\subsubsection{Iteration}
\label{c4}

Now classical arguments show the convergence of the iterative scheme
if $\eps \in ]0,\eps_0]$ and $\eps_0$ is small enough. We describe
the main lines (see \cite{S1}). Let us fix an integer $m
 > 4$, and note
$$
R := 1 + \sup_{0<\eps<1}  \{ \eps  \, ( \ \| w^{\eps, 0}\|_\infty
+\sum_{0\leq j\leq \mu} \|\cZ_j w^{\eps, 0}\|_\infty+\|\eps \D_x
w^{\eps, 0}\|_\infty \ ) \}.
$$

\begin{proposition}

Let be given $\lambda >1 $. 
Then there exists  $h>1$ such that  for  $\eps_0>0$ small enough,  for
all $\nu \in \NN$, for all $\eps \in ]0,\eps_0]$, there hold
\begin{equation}\label{53}
 \| w^{\eps, \nu}\|_\infty +\sum_{0\leq j\leq \mu}
\|\cZ_j w^{\eps, \nu}\|_\infty+\|\eps \D_x w^{\eps, \nu}\|_\infty  <
R\eps ^{-1}
\end{equation}
and
\begin{equation}\label{54}
 \| w^{\eps,\nu}\|_{m,\lambda} +
 \| \eps \D_\n w^{\eps,\nu}\|_{m,\lambda} <
h.
\end{equation}
\end{proposition}

\begin{proof} 

For  $h$ large enough,
the inequalities  (\ref{53}) and (\ref{54}) are satisfied for $\nu=0$. 
Now suppose that $w^{\eps,\nu}$ satisfies (\ref{53}), (\ref{54}). 
We want to prove that  $w^{\eps,\nu+1}$ also satisfies (\ref{53}), (\ref{54}).
The proposition \ref{22} gives a constant $C_m(R)$ and the inequality
(\ref{20}) holds with $\bu = w^{\eps,\nu+1}$, $\bb = w^{\eps,\nu}$,
and $f=f^{\eps,\nu}$ defined in (\ref{91}).
 In order to control the right hand side of (\ref{50}), we need a
control of $\|\cH(w^{\eps,\nu})\|_\infty$ and of
$\|\cH(w^{\eps,\nu})\|_{m,\lambda}$, which is a consequence of the
following lemma.
\begin{lemma}\label{94} Let $m\in \NN$. There exists $c>0$ such that for
all $\lambda \geq 1$,
\begin{equation}
\| \cH(v)\|_{m,\lambda} + \| \eps \D_\n \cH(v)\|_{m-1,\lambda} \leq c (\|
v\|_{m,\lambda} + \|\eps  \D_\n v\|_{m-1,\lambda} ).
\end{equation}
\end{lemma}

\begin{proof}
 We note $E(\D) := ( \div , \rot )$ the operator from $[\cS'(\RR^3)]^3$ to $[\cS'(\RR^3)]^4$. We denote
by $E^{-1}(\D)$ the inverse operator. Then $u = E^{-1}(\D) f$, is
defined by $\hat{u}(\xi) = -i |\xi|^{-2} \big(\hat{a}(\xi) \xi - \xi \wedge
\hat{b}(\xi) \big) $ where $\hat{f}(\xi)=( \hat{a} (\xi),
      \hat{b}(\xi)) \in \RR \times \RR^3$. Thus 
$\hat{u}(\xi) = M (\xi) \hat{f}(\xi)$,  where  $M(\xi)$ is a $3\times 4$ matrix whose entries are
\emph{rational functions of $\xi$ homogeneous of degree $-1$}. Let
us fix $\chi \in \cC^\infty_0(\RR^3, \RR)$ such that $\chi(\xi)=0$
when $|\xi| \leq 1$ and $\chi(\xi)=1$ when $|\xi|\geq 2$, and call
$P(D)$ and $R(D)$ the operators from  $[\cS'(\RR^3)]^4$ to
$[\cS'(\RR^3)]^3$ defined by
$P(D)f := \cF^{-1}\big(\, \chi M \hat{f}\, \big)$ and $ R(D)f :=
\cF^{-1}\big(\, (1-\chi) \hat{f}\, \big)$
where $\cF^{-1}$ means the inverse Fourier transform. In the
sequel we will simply note $\cS'(\RR^3)$ and $L^2(\tOmega)$ instead
of $\big[\cS'(\RR^3)\big]^4$ and $\big[L^2(\tOmega)\big]^4$, meaning
that we talk about the \emph{components} of the vector valued
functions, the (finite) number of components being understood. We
have $E^{-1}(\D) = P(D) + R(D).$
The operator $P(D)$ is a special case of classical
pseudo-differential operator of class $S^{1}_{-1,0}(\RR^3 \times
\RR^3)$, elliptic, and $R(D)$ is an infinitely smoothing operator of
class $S^{-\infty}_{1,0}(\RR^3 \times \RR^3)$.

Let us now take into account the $t$ coordinate. Let us note
$\tOmega= ]0,T[\times \Omega$,  $\tGamma = ]0,T[\times \Gamma$ and
$\tSigma = ]0,T[\times \Sigma$. We extend the actions of $P$ and $R$
to the spaces of functions or distributions which depend also on $t$
like $L^2(\tOmega)$ or $C\big( [0,T],\cS'(\RR^3)\big)$, by
considering $t$ as a parameter so that $Pu(t,x) := P(D)u(t,.)(x)$.
Let $v \in H^m_{co}(\tOmega; \RR^4)$ such that $\Dn v \in
H^{m-1}_{co}(\tOmega; \RR^4)$.   Then $\cH(v)= u_{|\tOmega}$ where
$u\in L^2([0,T]\times \RR^3)$ is defined by
$   E(\D)u = \overline{E(\D)v} +
    (v_{|\tGamma}.\bn) \otimes \delta_{\tGamma} $,
where the notation $\overline{V}$ means the extension of $V$ by $0$
to $[0,T]\times \RR^3$ (or to $\RR^3$, depending on the context).

Let us note
$f :=\overline{E(\D)v}$, which is in $H^{m-1}_{co}(]0,T[\times \RR^3),$
and $g=(v.\bn)_{|\Gamma}$. This trace is well defined since by
assumption $v\in H^1(\tOmega)$, and using local coordinates patches
one sees that
$g \in H^{m-\frac{1}{2}}(\tGamma)$,
the usual Sobolev spaces. The operator $P(D)$ satisfies the
\emph{transmission property} (introduced by Boutet de Monvel
\cite{BMV}, \cite{BMV2}) on $\Omega$ and on $\RR^3\setminus \Omega$
because its symbol is a rational function of $\xi$, which is a
sufficient condition to satisfy the transmission condition. The
transmission property has been also studied and used by  Grubb,
and we also refer to papers 
\cite{Grubb1} and \cite{Grubb2}. To avoid many repetitions, we will
note in what follows $\Omega_1 := \Omega$ and
$\Omega_2=\RR^3\setminus \Omega$. Since $P(D)$ is elliptic of order
1, the transmission property implies (see \cite{Grubb1} and
\cite{Grubb2}) that if $ v \in H^s(\Omega)$ then for $j=1,2$,
    $\big( P(D) \overline{v} \big)_{|\Omega_j}
    \in H^{s+1}(\Omega_j)$.
Let us note note $u_{(j)} = u_{|\tOmega_j}$, for $j=1, 2$, so that
$\cH(v)= u_{(1)}\in L^2(\tOmega)$. Using the notations of
\cite{BMV}, \cite{Grubb1}, \cite{Grubb2},
\begin{equation}\label{102}
u _{(j)}= \big( E^{-1}(D)\overline{v}\big)_{|\tOmega_{(j)}} =
P(D)^{(j)} f + K^{(j)}_\Gamma (g) + R(D)^{(j)}E(\D)\overline{v},
\end{equation}
where
 $P(D)^{(j)} f = (P(D)f)_{|\tOmega_{(j)}}$,
 $R(D)^{(j)}\ov= (R(D)\ov)_{|\tOmega_{(j)}}$ and where
$K^{(j)}_\Gamma (g) = \big(\,P(D)\big( g \otimes \delta_\Gamma)\,
\big)_{|\tOmega_{(j)}}$  is the "Poisson operator":
\begin{equation}\label{103}
K_\Gamma^{(j)} : H^s(\Gamma) \rightarrow H^{s+1/2}(\Omega^{(j)}),
\quad
\end{equation}
(linear continuous), extended to functions depending on $t$ as a
parameter. (See theorems 2.4 and 2.5 of \cite{Grubb2}).

Let us now prove the lemma. First of all, $\D_t^m\cH(v)=\cH(\D_t^m
v)$ is in $L^2(\tOmega)$ because $\D_t^m v\in L^2(\tOmega)$ and
$\cH$ acts on $L^2(\tOmega)$. It is also easy to show that
$\D_t^{m-1}\cH(v) \in H^1(\tOmega)$: by assumption, for any $t\in
[0,T]$, $\D_t^{m-1}v(t,.) \in H^1(\Omega)$, hence
$\cH(\D_t^{m-1}\ov)(t,.) \in H^1(\Omega)$ because
$\D_t^{m-1}\ov(t,.)$ is piecewise-$H^1$ and because of the
properties of $\cH$.
 Hence $\D_x\D_t^{m-1}\cH(v)
\in L^2(\tOmega)$ and since we already know that $\D_t^m\cH(v) \in
L^2(\tOmega)$ we have proved that $\D_t^{m-1}\cH(v)\in
H^1(\tOmega)$.

Let us show now that $\cZ_j\D_t^{m-2}\cH(v) \in H^1(\tOmega)$ for
$j=1, \dots , \mu$. We have
$    E(\D)u = f + g \otimes \delta_\Gamma $.
Since $E(\D)$ is elliptic (as an operator in $\cS'(\RR^3)$, but not
in $\cS'(\RR^4)$), we can express the normal derivatives of $u$ in
term of tangential derivatives and of $E(\D)u$, and this implies
that the commutator $[E(\D),\cZ_j]u$ writes
\begin{equation}\label{107}
[E(\D), \cZ_j ]u = \sum_1^\mu A_j \cZ_j u + A_0 f + B g\otimes
\delta_\Gamma
\end{equation}
where $A_j,B$ are matrices with $\cC_b^\infty$ entries (depending on
the fields $\cZ_j$). It follows that
$$
E(\D)\cZ_j u = \sum_{|\alpha|\leq 1} M_\alpha \cZ^\alpha f +
\sum_{|\alpha|\leq 1} N_\alpha (\cZ^\alpha g) \otimes \delta_\Gamma
$$
with $\cC_b^\infty(\RR^3)$ matrices $M_\alpha, N_\alpha$, and
applying $\D_t^{m-2}$ gives:
\begin{equation}\label{108}
E(\D)\cZ_j \D_t^{m-2} u = \sum_{|\alpha|\leq 1} M_\alpha .
\cZ^\alpha \D_t^{m-2}f + \sum_{|\alpha|\leq 1} N_\alpha .
(\cZ^\alpha\D_t^{m-2} g) \otimes \delta_\Gamma
\end{equation}
Now $\cZ^\alpha \D_t^{m-2}f \in L^2(\tOmega)$, because $f=
\overline{E(\D)v}$, and the transmission property implies that for every
$t\in [0,T]$, the function 
  $  P(D)^{(j)}\big(\, \cZ^\alpha\D_t^{m-1}f \, \big)(t,.)$
  is in $  H^1\big( \, \Omega_{(j)} \, \big)$. This implies that $\D_x P(D)^{(j)}\big(\,
\cZ^\alpha\D_t^{m-2}f \, \big)\in L^2(\tOmega)$ and since we already
know that $\D_t P(D)^{(j)}\big(\, \cZ^\alpha\D_t^{m-2}f \, \big)\in
L^2(\tOmega)$ from the previous case, we deduce that for $j=1, 2$ the functions 
    $P(D)^{(j)}\big(\, \cZ^\alpha\D_t^{m-2}f \, \big)$ is in
    $H^1\big( \, \tOmega_{(j)} \, \big)$.
Concerning the boundary term in (\ref{108}), since $g \in
H^{m-\frac{1}{2}}(\tGamma)$ we know that
 $\cZ^\alpha\D_t^{m-2} g \in H^{1/2}(\Gamma)$ and
the property (\ref{103}) implies that, for all $t\in [0,T]$, the
functions $K_\Gamma^{(j)} \big( \, \cZ^\alpha\D_t^{m-2} g \, \big)(t,.) $
is in
$H^1\big(\Omega_{(j)}\big)$.
By the same way as before we deduce that for $j=1, 2$ the functions 
$K_\Gamma^{(j)} \big( \, \cZ^\alpha\D_t^{m-2} g \, \big) $ is in
$H^1\big(\tOmega_{(j)}\big) $.
Now, applying $E(\D)^{-1} = P(D) + R(D)$ to the equation
(\ref{108}) gives $\cZ_j
\D_t^{m-2}u_{(j)} \in H^1\big(\tOmega_{(j)}\big)$ as claimed. Then,
the proof can be continued by induction in the same way.
\end{proof}

The lemma \ref{94},
together with the Gagliardo-Nirenberg-Moser estimates and the
induction assumption, implies that (like the majoration of the term
(5.25) in paper \cite{S1}):
\begin{equation}\label{95}
    \| f^{\eps,\nu} \|_{m,\lambda} \leq c (R)
    (\|w^{\eps,\nu}\|_{m,\lambda}+ \| \eps \D_x
    w^{\eps,\nu}\|_{m,\lambda}) < c(R) \rho(\lambda).
\end{equation}
Hence, the proposition \ref{22} implies that
\begin{equation}\label{96}
\begin{aligned}
\|\eps \D_x w^{\eps,\nu+1}\|_{m,\lambda}+ \lambda
\|w^{\eps,\nu+1}\|_{m,\lambda} &\leq C_m(R) \big[ \lambda^{-1} c(R)
\rho(\lambda) \\
&+ R (\|w^{\eps,\nu+1}\|_\infty + \| \eps
\D_xw^{\eps,\nu+1}\|_\infty)  + I_{m,\lambda}
(w^{\eps,\lambda})\big].
\end{aligned}
\end{equation}
We now use the following Sobolev inequalities:
\begin{eqnarray*}
    \eps^{1/2} \| \bu \|_\infty &\leq&  e^{\sigma \lambda }
    ( \| \bu\|_{m,\lambda} +
    \|\eps\Dn \bu\|_{m,\lambda}),
\\   \eps^{1/2} \| \eps\Dn \bu \|_\infty &\leq&  e^{\sigma \lambda }
    ( \| \eps \Dn\bu\|_{m,\lambda} +
    \|(\eps\Dn)^2\bu\|_{m,\lambda}).
\end{eqnarray*}
By taking $\lambda $ large enough et $\eps >0$
small enough 
the inequality (\ref{53}) is also satisfied for $w^{\eps, \nu +1}$
and the proof by induction is complete.
\end{proof}

Now by extracting a convergent subsequence it is a classical
argument to show the convergence in $L^2(]0,T[\times \Omega)$ of
$w^{\eps,\nu}$ to a solution $w^\eps$ of the non linear problem
which satisfies the same estimates $(\ref{53})$, $(\ref{54})$. This
concludes the proof of  Theorem $\ref{main2}$. \hfill $\Box$

\end{document}